\documentclass[5p]{elsarticle}

\usepackage[intlimits]{amsmath}
\usepackage{amsthm,amssymb}

\usepackage[english]{babel}
\usepackage[]{algorithm2e}
\SetKwProg{Init}{Initialization:}{}{}

\usepackage{graphics} 
\usepackage{epsfig} 
\usepackage{tikz}           
\theoremstyle{plain}
\newtheorem{assumption}{Assumption}

\newtheorem{lemma}{Lemma}
\newtheorem{theorem}{Theorem}
\newtheorem{proposition}{Proposition}
\theoremstyle{remark}

\biboptions{sort&compress}
\journal{Systems \& Control Letters}

\begin{document}

	\begin{frontmatter}
		
		\title{Observer Design for a Class of ODE - Continuum-PDE  Cascade Systems Inspired by a 
Control-Theoretic Model of Large-Scale Arterial Networks of Blood Flow\tnoteref{t1}\tnoteref{t2}}
		\tnotetext[t1]{Funded by the European Union (ERC, C-NORA, 101088147). Views and 
		opinions 
		expressed are however those of the authors only and do not necessarily reflect those of the 
		European Union or the European Research Council Executive Agency. Neither the European 
		Union nor the granting authority can be held responsible for them.}
		\tnotetext[t2]{MATLAB codes for the numerical simulations are available on 
		\texttt{https://github.com/jphumaloja/arterial-tree.}}
		
		\author[1]{Jukka-Pekka Humaloja}
		\author[1]{Nikolaos Bekiaris-Liberis}
		
		\affiliation[1]{organization={Department of Electrical and Computer Engineering, Technical 
		University of Crete},
			addressline={University Campus, Akrotiri}, 
			city={Chania},
			postcode={73100}, 
			country={Greece}}
		
\begin{abstract}
We develop a backstepping-based observer design for a class of ODE - continuum-PDE cascade 
systems, which can be viewed as the limit, of a finite collection of ODE - $2 \times 2$ hyperbolic 
systems, as the number of individual PDE system components tends to infinity. The large-scale 
collection of ODE - $2 \times 2$ hyperbolic systems is motivated by a dynamic model that we 
present, of a network of peripheral arteries, to which central (aortic) blood flow/pressure enters. 
We address a case in which average (boundary) measurements, over the ensemble dimension, 
are available, which is motivated by the availability of non-invasive, peripheral flow/pressure 
measurements. Exponential stability of the estimation error system is shown by proving 
well-posedness of the kernel equations and constructing a Lyapunov functional. We also establish 
that part of the backstepping kernels derived coincide with the solution of a Sylvester equation. 
We then 
apply the continuum-based observer for state estimation of the large-scale counterpart and, in 
particular, of the blood flow system, introducing an approach for optimal construction of 
continuum approximations. We also introduce an implementation method, adopting a 
spectral-based approach for computing the observer dynamics, which we illustrate in an 
academic, numerical simulation example. Furthermore, we illustrate the design in the problem of 
central flow/pressure estimation using realistic parameters and flow/pressure waveforms.
\end{abstract}
		
		
\end{frontmatter}

\section{Introduction}

Hyperbolic partial differential equation (PDE) continua can be 
either viewed as limits of 
large-scale hyperbolic PDEs, as the number of PDE subsystems tends to infinity \cite{HumBek25, 
HumBek25b, HumBek26,HumBek26c}, which may describe the dynamics of, for example, blood 
flow in parallel 
arterial networks \cite{SwaXuD09}, or traffic \cite{TumCan22}, \cite{MalPap21b} and water 
\cite{BasCorBook, JhaLel22} flows in large-scale networks; or may themselves serve as models 
for describing the effect of continuum parameters, such as, for example, drivers' age in 
multi/continuum-class traffic \cite{AllKrs25}. In particular, adopting a continuum PDE viewpoint 
enables computation of observer/control kernels and observer dynamics independently of the 
number of PDE 
subsystems, thus exhibiting the potential of reducing computational complexity of 
control/observer kernels and observer dynamics for large-scale systems. Motivated more 
specifically by the problem of 
aortic flow/pressure estimation employing blood flow measurements from peripheral arteries (as 
in 
\cite{SwaXuD09}), in the present paper we develop an observer design for a class of linear, 
ordinary differential equation (ODE) - continuum-PDE cascades; and we apply it for state 
estimation of the large-scale counterpart and, in particular, of the blood flow system.

Backstepping-based control and estimation of cascade systems involving ODEs and hyperbolic 
PDEs have been addressed in, for example, \cite{AndVaz18, AnfAam18, AurBri23, BekKrs10, 
BekKrs14, DeuGeh18, DeuGeh19, DiMArg18, DiMLam20, HasAam16, IrsEspIFAC22, IrsDeu23, 
RedAur21, 
VazAur26, WanKrsBook, XuXLiu24}. Alternative methods addressing the 
problem of control/observer design for ODE - hyperbolic PDE cascades can be found in, for 
example, \cite{TriAnd18, TerFri19, FerCri20, MarBri21, Nat21, MatWuY22, SelFri24, Lha26}. 
Because as the number of PDE 
subsystems becomes 
large, solving 
the respective backstepping kernel equations and computing the observer states may become 
computationally laborious, 
\cite{HumBek25, HumBek25b, HumBek26, HumBek26b,HumBek26c} introduced an approach in 
which 
backstepping 
kernels and observers can be computed independently of the number of PDE subsystems. This 
remedies the 
increase in computational complexity with the number of PDE subsystems. The approach relies on 
approximating the exact kernels/observers by respective continuum kernels/observers, as the 
kernels developed in 
\cite{AllKrs25} for a class of hyperbolic PDE continua. There is no result addressing the problems 
of observer design for linear ODE - hyperbolic PDE-continua and its application to state 
estimation of large-scale collections of ODE - $2\times 2$ hyperbolic PDEs. Furthermore, our 
modeling 
approach for blood flow in arterial networks is inspired by \cite{SwaXuD09}, while the model is 
derived utilizing \cite{SinBek24, StePhD, Bek23, ReyMer09} as theoretical basis. To the best of 
our knowledge, there is no control-theoretic modeling and estimation approach for the class of 
blood flow arterial networks considered here.

In the present paper, we consider a class of ODE - continuum-PDE cascade systems, for which  
we design a Luenberger-type continuum observer, with output injection gains designed based on 
a continuum-PDE backstepping transformation that we introduce. We show well-posedness of the 
respective backstepping kernels, employing pointwise in the ensemble variable arguments, which 
enables us to utilize the respective results for scalar $2\times 2$ hyperbolic systems from 
\cite{CorVaz13}. Because only average (over the ensemble dimension) boundary measurements 
are available, we derive a sufficient condition (that also depends on the choice of the target 
system as in, e.g., \cite{BasCorBook, BouBri20}) for exponential stability of the estimation error 
dynamics via construction of a suitable Lyapunov functional. As an alternative approach of 
constructing the observer kernels, we establish that the proposed continuum observer coincides 
with the 
Sylvester equation-based design from \cite{Nat21}, which has not been shown before for any 
ODE-PDE cascade system. We then adopt a spectral-based approach to implement the observer 
numerically, which preserves the continuum aspect of the PDE and does not resort to explicit 
discretization of the spatial domain, thus providing a computationally tractable implementation 
method. The proposed implementation is illustrated in an academic, 
numerical simulation example.

We then utilize the continuum-based observer design for state estimation of a large-scale system 
counterpart. In particular, we introduce an approach in which a continuum approximation of a 
given large-scale system is constructed in an optimal manner; namely, such that the difference 
between the parameters of the continuum system and the parameters of the large-scale 
counterpart is minimized in $L^2$ sense. We establish that the dynamics of the estimation error 
of the state of the large-scale system under the continuum-based observer, satisfies an 
output-to-state stability-type estimate with a gain that depends on the continuum approximation 
error of the parameters. We specialize the design approach to aortic flow/pressure estimation 
using non-invasive measurements from peripheral arteries. We provide consistent simulation 
results for both the cases of central flow and pressure estimation, using realistic parameters and 
flow/pressure waveforms taken from \cite{StePhD}.

A preliminary version of this paper has been accepted to 2026 IEEE Conference on Decision 
and Control \cite{HumBekCDC26}, 
containing preliminary versions of Sections~\ref{sec:obs} and~\ref{num:spec}, and results that are 
restricted to  continuum observer design for ODE - continuum-PDE cascades and its
numerical 
implementation. In the present paper, we extend these preliminary results in that, here, we: 
\begin{itemize}
\item[a)] utilize 
the continuum-based observer design for state estimation of large-scale system counterparts and 
derive an output-to-state stability-type estimate for the estimation error;
\item[b)] introduce an approach 
for construction of a continuum approximation for a given large-scale system, which is optimal in 
the $L^2$ sense; and
\item[c)] apply the developed design approach to aortic flow and pressure estimation, 
including presentation of consistent simulation results for both cases, using realistic parameter 
and waveform data. 
\end{itemize}

The rest of the paper is organized as follows. In Section~\ref{sec:obs}, we present the class of 
ODE - continuum-PDE cascade systems considered, design an observer, show stability of the 
estimation 
error dynamics, and establish the connection to the Sylvester equation-based design. In 
Section~\ref{sec:art}, we
introduce a 
constructive  approach to approximate ODE - large-scale $2\times2$ hyperbolic PDE 
cascades by ODE - 
continuum-PDE cascades
and show that the respective continuum observer can (approximately) estimate the state of 
the ODE - large-scale PDE cascade. In Section~\ref{sec:num}, we present a spectral-based 
implementation of the continuum observer that we illustrate in numerical simulations. In 
Section~\ref{sec:artappl}, we present a model of a large-scale, blood flow arterial 
network and apply the proposed continuum observer to aortic flow and pressure estimation.
Section~\ref{sec:conc} concludes the paper.

\paragraph*{Notation} 
We denote by
$E_c = L^2([0,1]; L^2([0,1]; \mathbb{R}))$ a continuum space  equipped with the 
inner product
\begin{align}
\langle f_1, f_2 \rangle_{E_c} = \int\limits_0^1\int\limits_0^1 f_1(x,\zeta)f_2(x,\zeta)d\zeta dx.
\end{align}
Hence, $E_c^2$ can be viewed as the continuum limit of the 
space $E = L^2([0,1]; \mathbb{R}^{2m})$ equipped with the inner product
\begin{align}
\left\langle \left(\begin{smallmatrix}
\mathbf{u}_1 \\ \mathbf{v}_1
\end{smallmatrix}\right), \left(\begin{smallmatrix}
\mathbf{u}_2 \\ \mathbf{v}_2
\end{smallmatrix}\right) \right\rangle_E & = \nonumber \\
\int\limits_0^1\frac{1}{m} \sum_{i=1}^m u_1^i(x)u_2^i(x)dx + \int\limits_0^1\frac{1}{m} 
\sum_{i=1}^m v_1^i(x)v_2^i(x)dx,
\end{align}
for some $m \in \mathbb{N}$, as $m \to \infty$. The transform 
$\mathcal{F}_m$ maps any $\mathbf{b} \in \mathbb{R}^m$ to $L^2([0,1]; \mathbb{R})$ as
\begin{equation}
	\label{eq:Fm}
	\mathcal{F}_m\mathbf{b} =\sum_{i=1}^m b_i \chi_{((i-1)/m,i/m]},
\end{equation}
where $\chi_{((i-1)/m,i/m]}$ denotes the indicator function of the interval $((i-1)/m,i/m]$. 
Moreover, we denote by $\mathcal{T}$ the 
triangular set
\begin{align}
	\mathcal{T} & = \left\{ (x,\xi) \in [0,1]^2: \xi \leq x\right\}.
\end{align}
Finally, we say that a 
system is exponentially stable on a normed space $Z$ if, for any initial condition $z_0\in Z$, the 
(weak) solution 
$z(t)$ of the system satisfies $\|z(t)\|_Z \leq Me^{-\omega t}\|z_0\|_E$ for some $M,\omega > 0$ 
that are independent of $z_0$.

\section{Observer Design for a Class of ODE - Continuum-PDE Cascade Systems} 
\label{sec:obs}

\subsection{ODE - Continuum-PDE Cascade Systems}

The considered class of ODE - Continuum-PDE cascade systems is of the form
\begin{subequations}
\label{eq:infart}%
\begin{align}
\dot{\mathbf{X}}(t) & = \mathbf{A}\mathbf{X}(t), \\
u_t(t,x,y) + \lambda(x,y)u_x(t,x,y)   & = W(x,y)v(t,x,y), \label{eq:infloc1} \\
v_t(t,x,y) - \mu(x,y)v_x(t,x,y)
& = \theta(x,y)u(t,x,y), \label{eq:infloc2}
\end{align}
\end{subequations}
with boundary conditions 
\begin{subequations}
\label{eq:infartbc}%
\begin{align}
u(t,0,y) & = Q(y)v(t,0,y), \\
v(t,1,y) & = R(y)u(t,1,y) + F(y)\mathbf{C}\mathbf{X}(t),
	\end{align}
\end{subequations}
and a weighted averaged measurement 
\begin{equation}
\label{eq:out}
Y(t) = \int\limits_{y_1}^{y_2}g(y)v(t,0,y)dy,
\end{equation}
over some $0 \leq y_1 < y_2 \leq 1$ with weight $g$. 

We note that $y$ is viewed as an ensemble (index) variable rather than as a spatial 
coordinate (and no spatial derivatives in $y$ are involved). Hence, 
the PDE in \eqref{eq:infart} can be viewed as quasi-2-D or 1.5-D, as also discussed in \cite[Sect. 
I.B]{AllKrs25}. Regardless, the system \eqref{eq:infart}, \eqref{eq:infartbc} is interpreted in the 
weak sense on the state space $\mathbb{R}^n \times E_c^2$.

The parameters of \eqref{eq:infart}--\eqref{eq:out} satisfy the following assumption.
\begin{assumption}
\label{ass:loc}
In \eqref{eq:infart}, \eqref{eq:infartbc}, $\mathbf{A} \in \mathbb{R}^{n\times n}$ and $\mathbf{C} 
\in \mathbb{R}^{1\times n}$, where the eigenvalues of $\mathbf{A}$ are of geometric multiplicity 
one\footnote{This assumption is satisfied by the types of ODEs considered here, but it can be 
omitted by employing generalized eigenvectors.}. The PDE parameters satisfy $\lambda,\mu \in 
C^1([0,1]^2; \mathbb{R})$, $w, \theta, \in C([0,1]^2; \mathbb{R})$, and $Q,R,F\in C([0,1]; 
\mathbb{R})$, with $\lambda(x,y)$, $\mu(x,y) > 0$ for all $x,y \in [0,1]$ and $|Q(y)R(y)| < 1$ for all 
$y \in [0,1]$. In \eqref{eq:out}, $g \in L^2([0,1]; \mathbb{R})$\footnote{For modeling a point 
measurement, the choice
$g(y) = \delta(y-y_0)$ for some $y_0 \in [0,1]$ is also possible, where $\delta$ denotes the Dirac 
delta function.}.
\end{assumption}
\noindent In order for the system \eqref{eq:infart}, \eqref{eq:infartbc} to be observable, the signal 
$\mathbf{CX}(t)$ has to be sufficiently ``visible'' in the measurement $Y(t)$. In
technical terms, there cannot be transmission zeros from $\mathbf{CX}(t)$ to $Y(t)$, which can 
be formulated as follows (see, e.g., \cite[Assum. 4.2, Prop. 4.7]{Nat21}).
\begin{assumption}
\label{ass:obs}
Denote by $G$ the transfer function from $\mathbf{CX}$ to $Y$. Assume that $G(s)$ exists for 
all $s \in \sigma(\mathbf{A})$, while $G(s)\mathbf{C}\mathbf{v} \neq 0$ for all $s \in 
\sigma(\mathbf{A})$ and the corresponding eigenvector $\mathbf{v}$ satisfying $\mathbf{Av} = 
s\mathbf{v}$.
\end{assumption}

If detectability was considered instead of observability, Assumption~\ref{ass:obs} would only have 
to concern the unstable eigenvalues of $\mathbf{A}$ and the respective eigenvectors. Note that 
this only applies to the ODE part of \eqref{eq:infart}, \eqref{eq:infartbc}, so that the PDE part has 
to be observable/detectable based on the output $Y(t)$ as well, in order for the whole system to 
be observable/detectable. With full boundary measurements $v(t,0,y)$ (or $u(t,1,y)$), it follows 
adapting the arguments for a finite number of systems, 
e.g., from \cite[Sect. 8.3]{AnfAamBook}, that the PDE part is detectable for all PDE parameters 
satisfying Assumption~\ref{ass:loc}. However, additional assumptions are required if additional 
backstepping kernels are employed (compared to those in \cite[Sect. 8.3]{AnfAamBook}), e.g., as 
in \cite[Sect. IV.B]{VazKrsCDC11}, where $q \neq 0$ (translating to $Q(y) \neq 0$ for all $y \in 
[0,1]$ in \eqref{eq:infartbc}) is required for well-posedness of the additional kernels. 

In our case, the well-posedness of the employed backstepping kernels (for the chosen target 
system) requires that $R(y) \neq 0$ 
for all $y \in [0,1]$ that can be viewed as an observability assumption of the PDE part when the 
full, boundary measurement $v(t,0,y)$ is available (see, e.g., \cite{LiTRao10}), albeit 
the implementation of the proposed observer is independent of this assumption. Furthermore, 
since here the available measurement is only an average (over $y$) of the full, boundary 
measurement $v(t,0,y)$, additional assumptions are required to be imposed on the parameters of 
the PDE part to guarantee stability of the estimation error system (see Lemma~\ref{lem:stab}). 
These conditions may in fact imply open-loop stability of the PDE part of system \eqref{eq:infart}, 
\eqref{eq:infartbc}, 
however, it is not clear how they could be removed or alleviated in the present case in which only 
an average (over $y$) boundary measurement is available. Nevertheless, the main idea of the 
backstepping procedure we present would be the same even if a full, boundary measurement is 
available, in which case such a restrictive assumption could be removed; while as regards the 
application to the arterial network, open-loop stability (of each arterial component) is supported, 
in principle, by the results in \cite{SwaXuD09, SinBek24, StrPetCPDE25}. Finally, if the PDE 
part of system \eqref{eq:infart}, \eqref{eq:infartbc} is known to be stable, then the observer can 
be designed directly using the Sylvester equation-based approach from \cite[Sect. 
4]{Nat21} without considering a backstepping transformation. The Sylvester equation-based 
design, in fact, coincides with the backstepping one as shown in Section~\ref{obs:syl}.

\subsection{Continuum Observer Design}

We consider an observer design based on weighted averaged measurement \eqref{eq:out}. A 
typical Luenberger observer is of the form 
\begin{subequations}
\label{eq:lobs}
\begin{align}
\dot{\hat{\mathbf{X}}}(t) & = \nonumber \\
\mathbf{A}\hat{\mathbf{X}}(t) + \mathbf{L}\left(\int\limits_{y_1}^{y_2}g(y)\hat{v}(t,0,y)dy - 
Y(t)\right), \\
\hat{u}_t(t,x,y) + \lambda(x,y)\hat{u}_x(t,x,y) - W(x,y)\hat{v}(t,x,y)& = \nonumber \\
 P_1(x,y)\left(\int\limits_{y_1}^{y_2}g(y)\hat{v}(t,0,y)dy - Y(t)\right), 
\label{eq:lobspde1}
\end{align}
\begin{align}
\hat{v}_t(t,x,y) - \mu(x,y)\hat{v}_x(t,x,y) - \theta(x,y) \hat{u}(t,x,y) & = \nonumber \\
P_2(x,y)\left(\int\limits_{y_1}^{y_2}g(y)\hat{v}(t,0,y)dy - Y(t)\right), 
\label{eq:lobspde2}
\end{align}
\end{subequations}
with boundary conditions
\begin{subequations}
\label{eq:lobsbc}
\begin{align}
\hat{u}(t,0,y) & = Q(y)\hat{v}(t,0,y), \label{eq:lobsbc1} \\
\hat{v}(t,1,y) & = R(y)\hat{u}(t,1,y) + F(y)\mathbf{C}\hat{\mathbf{X}}(t).
\end{align}
\end{subequations}
The gain $\mathbf{L} \in \mathbb{R}^n$ is chosen such that the matrix $\mathbf{A} + 
\mathbf{L}\int\limits_{y_1}^{y_2}g(y)\pmb{\gamma}_0(y)dy$ is Hurwitz (which is always possible 
under Assumption~\ref{ass:obs}), and $P_1,P_2 \in 
C([0,1]^2; \mathbb{R})$ are chosen as
\begin{subequations}
\label{eq:P12}
\begin{align}
P_1(x,y) & = \pmb{\gamma}^1(x,y)\mathbf{L},  \\
P_2(x,y) & = \pmb{\gamma}^2(x,y)\mathbf{L},
\end{align}
\end{subequations}
where $\pmb{\gamma}_0 \in C([0,1]; \mathbb{R}^{1\times n})$ is given by
\begin{align}
\label{eq:g0}
 \pmb{\gamma}_0(y) & = F(y)\mathbf{C}\left(\begin{bmatrix} Q(y)I_{n\times n} & I_{n\times n} 
\end{bmatrix} \right. \nonumber \\
& \qquad \times \exp\left(\int\limits_0^1\begin{bmatrix}
-\frac{\mathbf{A}}{\lambda(\xi,y)} & -\frac{\theta(\xi,y)}{\mu(\xi,y)}I_{n\times n} \\
\frac{W(\xi,y)}{\lambda(\xi,y)}I_{n\times n}  & \frac{\mathbf{A}}{\mu(\xi,y)}
\end{bmatrix}d\xi\right) \nonumber \\
& \qquad \left. \times \begin{bmatrix}
-R(y)I_{n\times n} \\ I_{n\times n}
\end{bmatrix}\right)^{-1},
\end{align}
and $(\pmb{\gamma}^1, \pmb{\gamma}^2) =: \pmb{\gamma} \in C([0,1]^2; \mathbb{R}^{1\times 
2n})$ is given by
\begin{align}
\label{eq:gamma}
\pmb{\gamma}(x,y) & = \begin{bmatrix}
Q(y)\pmb{\gamma}_0(y) & \pmb{\gamma}_0(y)
\end{bmatrix} \nonumber \\
& \quad \times \exp\left(\int\limits_0^x\begin{bmatrix}
-\frac{\mathbf{A}}{\lambda(\xi,y)} & -\frac{\theta(\xi,y)}{\mu(\xi,y)}I_{n\times n} \\
\frac{W(\xi,y)}{\lambda(\xi,y)}I_{n\times n} & \frac{\mathbf{A}}{\mu(\xi,y)}
\end{bmatrix}d\xi\right).
\end{align}
We note that the inverse appearing in \eqref{eq:g0} is well-defined under 
Assumptions~\ref{ass:loc} and~\ref{ass:obs}, as shown in the proofs of Lemma~\ref{lem:bswp} 
and Proposition~\ref{prop:sylv}.

\subsection{Stability of Continuum Estimation Error System}

We start presenting the estimation error dynamics for
$\tilde{\mathbf{X}} = \hat{\mathbf{X}}-\mathbf{X}, \tilde{u} = \hat{u}-u, \tilde{v}=\hat{v}-v$ as
\begin{subequations}
\label{eq:lerr}
\begin{align}
\dot{\tilde{\mathbf{X}}}(t) & = \nonumber \\
\mathbf{A}\tilde{\mathbf{X}}(t) + \mathbf{L}\int\limits_{y_1}^{y_2}g(y)\tilde{v}(t,0,y)dy, \\
\tilde{u}_t(t,x,y) + \lambda(x,y)\tilde{u}_x(t,x,y) & = \nonumber \\
W(x,y)\tilde{v}(t,x,y) + P_1(x,y)\int\limits_{y_1}^{y_2}g(y)\tilde{v}(t,0,y)dy, 
 \\
\tilde{v}_t(t,x,y) - \mu(x,y)\tilde{v}_x(t,x,y) & = \nonumber \\
\theta(x,y) \tilde{u}(t,x,y) + P_2(x,y)\int\limits_{y_1}^{y_2}g(y)\tilde{v}(t,0,y)dy, 
\end{align}
\end{subequations}
with boundary conditions
\begin{subequations}
\label{eq:lerrbc}
\begin{align}
\tilde{u}(t,0,y) & = Q(y)\tilde{v}(t,0,y),
\label{eq:lerrbc1} \\
\tilde{v}(t,1,y) & = R(y)\tilde{u}(t,1,y) + F(y)\mathbf{C}\tilde{\mathbf{X}}(t).
\end{align}
\end{subequations}

Consider a backstepping transformation of the form
\begin{subequations}
  \label{eq:obsbs}
\begin{align}
  \tilde{u}(t,x,y)
  & = \alpha(t,x,y) - \int\limits_0^x
                 N^{11}(x,\xi,y)\alpha(t,\xi,y)d\xi \nonumber \\
  & \qquad - \int\limits_0^xN^{12}(x,\xi,y)\beta(t,\xi,y)d\xi -
    \pmb{\gamma}^1(x,y)\tilde{\mathbf{X}}(t), \\
  \tilde{v}(t,x,y)
  & = \beta(t,x,y) -
    \int\limits_0^xN^{21}(x,\xi,y)\alpha(t,\xi,y)d\xi \nonumber \\
 & \qquad  - \int\limits_0^xN^{22}(x,\xi,y)\beta(t,\xi,y)d\xi - 
 \pmb{\gamma}^2(x,y)\tilde{\mathbf{X}}(t),
\end{align}
\end{subequations}
where $N^{11}, N^{12}, N^{21}, N^{22} \in C(\mathcal{T} \times [0,1]; \mathbb{R})$ are yet to be 
determined, $\pmb{\gamma}^1, \pmb{\gamma}^2$ are given by \eqref{eq:g0}, \eqref{eq:gamma}, 
and $(\tilde{\mathbf{X}}, \alpha, \beta)$ denotes the estimation error target system with dynamics
\begin{subequations}
\label{eq:obstsmod}
\begin{align}
\dot{\tilde{\mathbf{X}}}(t) & = 
\left(\mathbf{A} + \mathbf{L}\int\limits_{y_1}^{y_2}g(y)\pmb{\gamma}^2(0,y)dy 
\right)\tilde{\mathbf{X}}(t) \nonumber \\
& \qquad -
\mathbf{L}\int\limits_{y_1}^{y_2}g(y)\beta(t,0,y)dy, \\
\alpha_t(t,x,y) & = -\lambda(x,y)\alpha_x(t,x,y) + G_1(x,y)\beta(t,0,y), \\
\beta_t(t,x,y) & =  \mu(x,y)\beta_x(t,x,y) + G_2(x,y)\beta(t,0,y),
\end{align}
\end{subequations}
with boundary conditions 
\begin{subequations}
\label{eq:obstsbc}
\begin{align}
\alpha(t,0,y) & = Q(y)\beta(t,0,y), \\
\beta(t,1,y) & = R(y)\alpha(t,1,y),
\end{align}
\end{subequations}
where
\begin{subequations}
\label{eq:G12}
\begin{align}
G_1(x,y) & = N^{11}(x,0,y)\lambda(0,y)Q(y) - N^{12}(x,0,y)\mu(0,y), \\
G_2(x,y) & = N^{21}(x,0,y)\lambda(0,y)Q(y) - N^{22}(x,0,y)\mu(0,y).
\end{align}
\end{subequations}

Stability of the estimation error system \eqref{eq:lerr} follows directly from the following two 
lemmas.
\begin{lemma}[Well-posedness of backstepping kernels]
\label{lem:bswp}
Under Assumptions~\ref{ass:loc} and \ref{ass:obs}, if $R(y) \neq 0$ for all $y \in [0,1]$, then  
the backstepping transformation 
\eqref{eq:obsbs} mapping \eqref{eq:lerr}, \eqref{eq:lerrbc} to 
\eqref{eq:obstsmod}, \eqref{eq:obstsbc} is well-defined. In particular, $N^{11}, N^{12}, N^{21}$, 
$N^{22} \in C(\mathcal{T} \times [0,1]; \mathbb{R})$ and $\pmb{\gamma}_1, \pmb{\gamma}_2 \in 
C([0,1]^2; \mathbb{R}^{1\times n})$.

\begin{proof}
Setting \eqref{eq:lerr}, \eqref{eq:lerrbc} and \eqref{eq:obstsmod}, \eqref{eq:obstsbc} equal, we 
get that the kernels need to satisfy (see \ref{app:ker} for detailed derivations)
\begin{subequations}
\label{eq:obsk}
\begin{align}
\lambda(x,y)N_x^{11}(x,\xi,y) + N_\xi^{11}(x,\xi,y)\lambda(\xi,y) & = \nonumber \\
 -\lambda_\xi(\xi,y)N^{11}(x,\xi,y) + W(x,y)N^{21}(x,\xi,y), \\
\lambda(x,y)N_x^{12}(x,\xi,y) - N_\xi^{12}(x,\xi,y)\mu(\xi,y) & = \nonumber \\
N^{12}(x,\xi,y)\mu_\xi(\xi,y) + W(x,y)N^{22}(x,\xi,y), \\
\mu(x,y)N_x^{21}(x,\xi,y) - N_\xi^{21}(x,\xi,y)\lambda(\xi,y) & = \nonumber \\
N^{21}(x,\xi,y)\lambda_\xi(\xi,y) - \theta(x,y)N^{11}(x,\xi,y), \\
\mu(x,y)N_x^{22}(x,\xi,y) + N_\xi^{22}(x,\xi,y)\mu(\xi,y) & = \nonumber \\
-N^{22}(x,\xi,y)\mu_\xi(\xi,y) - \theta(x,y)N^{12}(x,\xi,y), \\
\lambda(x,y)\pmb{\gamma}_x^1(x,y) + \pmb{\gamma}^1(x,y)\mathbf{A} - 
W(x,y)\pmb{\gamma}^2(x,y)
 & = 0, \label{eq:kg1} \\
\mu(x,y)\pmb{\gamma}_x^2(x,y) - \pmb{\gamma}^2(x,y)\mathbf{A} + 
\theta(x,y)\pmb{\gamma}^1(x,y) 
 & = 0, \label{eq:kg2}
\end{align}
\end{subequations}
with boundary conditions
\begin{subequations}
\label{eq:obskbc}
\begin{align}
N^{12}(x,x,y)\mu(x,y) + \lambda(x,y)N^{12}(x,x,y) & = -W(x,y), \\
N^{21}(x,x,y)\lambda(x,y) + \mu(x,y)N^{21}(x,x,y) & = \theta(x,y), \\
N^{21}(1,\xi,y) - R(y)N^{11}(1,\xi,y) & = 0, \label{eq:obskbcN11} \\
N^{22}(1,\xi,y)  - R(y)N^{12}(1,\xi,y) & = 0, \label{eq:obskbcN12} \\
Q(y)\pmb{\gamma}^2(0, y) & = 
\pmb{\gamma}^1(0,y), \label{eq:kg1bc} \\
R(y)\pmb{\gamma}^1(1,y) + F(y)\mathbf{C} & =  \pmb{\gamma}^2(1,y). \label{eq:kg2bc}
\end{align}
\end{subequations}
Note that the equations \eqref{eq:kg1}, \eqref{eq:kg2} and boundary conditions \eqref{eq:kg1bc}, 
\eqref{eq:kg2bc} for $\pmb{\gamma}^1, \pmb{\gamma}^2$ are decoupled from $N^{11}, N^{12}$, 
$N^{21}, N^{22}$. Moreover, direct computation shows that the solution to  \eqref{eq:kg1}, 
\eqref{eq:kg2} with \eqref{eq:kg1bc}, \eqref{eq:kg2bc} is given by \eqref{eq:g0}, 
\eqref{eq:gamma}, which is well-defined by (the first part of) Assumption~\ref{ass:obs} which 
guarantees that the matrix inverse in \eqref{eq:g0} is well-defined\footnote{This is more apparent 
based on the alternative, spectral form \eqref{eq:gammas}--\eqref{eq:vk}.}. Finally, due to 
continuity of the parameters, the well-posedness of $N^{11}, N^{12}, N^{21}, 
N^{22}$ kernels can be argued pointwise in $y \in [0,1]$. That is, fixing arbitrary $y_0 \in [0,1]$ 
and introducing auxiliary kernels $\bar{N}^{ij}(\chi,\zeta) = 
N^{ij}(1-\zeta,1-\chi,y_0) = N^{ij}(x,\xi,y_0)$, for $i,j \in \{1,2\}$, the well-posedness of 
$\bar{N}^{ij}$ follows from \cite[Thm A.1]{CorVaz13} for parameters evaluated at $y = y_0$. 
Hence, the kernels $N^{ij}$ are well-posed pointwise in $y$. Moreover, the kernels at $y = y_0$ 
depend continuously on the parameter values at $y = y_0$, which together with the continuity of 
the parameters in $y$ implies that the kernels $N^{ij}$ are continuous in $y$.
\end{proof}
\end{lemma}

\begin{lemma}[Stability of target system] \label{lem:stab}
Under the conditions of Lemma~\ref{lem:bswp}, if 
\begin{equation}
M_G = \max\left\{\max_{x,y \in 
[0,1]}|G_1(x,y)|,\max_{x,y \in [0,1]}|G_2(x,y)|\right\},
\end{equation}
satisfies
\begin{equation}
\label{eq:MGcond}
\resizebox{.98\columnwidth}{!}{$\displaystyle 
M_G < \min_{y \in [0,1]}\sqrt{\frac{1}{2}\frac{1 - 
M_{RQ}^2\exp\left(\delta\int\limits_0^1\left(\frac{1}{\lambda(s,y)}
 + \frac{1}{\mu(s,y)}\right)ds\right)}{\frac{\exp\left(
 \frac{\delta}{m_\mu}\right)}{\delta m_\mu} + \frac{M_R^2}{\delta 
 m_\lambda}\exp\left(\delta\int\limits_0^1\left(\frac{1}{\lambda(s,y)}
 + \frac{1}{\mu(s,y)}\right)ds\right)}}$},
\end{equation}
where
\begin{subequations}
\begin{align}
M_{RQ} & = \max_{y \in [0,1]} |R(y)Q(y)|, & M_R & = \max_{y \in [0,1]} |R(y)|, \\
m_\lambda & = \min_{x,y \in [0,1]}\lambda(x,y), & m_\mu & = \min_{x,y \in [0,1]}\mu(x,y), 
\end{align}
\end{subequations}
and $\delta > 0$ is sufficiently small such that 
\begin{equation}
\label{eq:deltacond}
M_{RQ}^2\max_{y \in [0,1]}\exp\left(\delta\int\limits_0^1\left(\frac{1}{\lambda(s,y)}
 + \frac{1}{\mu(s,y)}\right)ds\right)  < 1,
\end{equation}
then, the target system \eqref{eq:obstsmod}, 
\eqref{eq:obstsbc} is exponentially stable on $\mathbb{R}^n \times E_c^2$.

\begin{proof}
Motivated by the proof of \cite[Thm 2.4]{BasCorBook} for scalar 
$2\times 2$ systems, we consider a Lyapunov functional for the PDE part of the form
\begin{align}
V(t) & = \int\limits_0^1\int\limits_0^1 
\frac{p_1(y)}{\lambda(x,y)}\alpha^2(t,x,y)\exp\left(-\delta\int\limits_0^x\frac{ds}{\lambda(s,y)}\right)dydx
\nonumber 
\\
& \quad + \int\limits_0^1\int\limits_0^1 
\frac{1}{\mu(x,y)}\beta^2(t,x,y)\exp\left(\delta\int\limits_0^x\frac{ds}{\mu(s,y)}\right)dydx,
\end{align}
for some positive (continuous) weight $p_1 > 0$ and constant $\delta > 0$. Thereafter, the total 
Lyapunov functional for \eqref{eq:obstsmod}, \eqref{eq:obstsbc} is taken of the form
\begin{equation}
\label{eq:Vtot}
V_\mathrm{tot}(t) = V(t) + \gamma \tilde{\mathbf{X}}^T(t)\mathbf{P}\tilde{\mathbf{X}}(t),
\end{equation}
for some $\gamma > 0$ and $\mathbf{P} > 0$ such that $\mathbf{A}_L^T\mathbf{P} + 
\mathbf{P}\mathbf{A}_L = -I_{n\times n}$,
where we denote $\mathbf{A}_L = \mathbf{A} + 
\mathbf{L}\int\limits_{y_1}^{y_2}g(y)\pmb{\gamma}^2(0,y)dy$, which follows under 
Assumption~\ref{ass:obs} since from \cite[Prop. 4.7]{Nat21} the pair $\left(\mathbf{A}, 
\int\limits_{y_1}^{y_2}g(y)\pmb{\gamma}^2(0,y)dy\right)$ is detectable. Computing 
the derivative of $V_\mathrm{tot}$ along the solution to \eqref{eq:obstsmod}, 
\eqref{eq:obstsbc} gives
\begin{align}
\dot{V}_{\mathrm{tot}}(t) & = -\delta V(t) - \int\limits_0^1(1- p_1(y)Q^2(y))\beta^2(t,0,y)dy 
\nonumber \\
& \quad -\int\limits_0^1 \left(p_1(y)\exp\left(-\delta \int\limits_0^1\frac{ds
}{\lambda(s,y)}\right) \right. \nonumber \\
& \left.\qquad \qquad - 
R^2(y)\exp\left(\delta\int\limits_0^1\frac{ds}{\mu(s,y)}\right)\right)\alpha^2(t,1,y)dy 
\nonumber \\
& \quad +\resizebox{.79\columnwidth}{!}{$\displaystyle \int\limits_0^1\int\limits_0^1 
\frac{2p_1(y)}{\lambda(x,y)}\alpha(t,x,y)G_1(x,y)\beta(t,0,y)\exp\left(-\delta\int\limits_0^x\frac{ds 
}{\lambda(s,y)}\right)dy dx$} \nonumber \\
& \quad +\resizebox{.79\columnwidth}{!}{$\displaystyle\int\limits_0^1\int\limits_0^1 
\frac{2}{\mu(x,y)}\beta(t,x,y)G_2(x,y)\beta(t,0,y)\exp\left(\delta\int\limits_0^x\frac{ds 
}{\mu(s,y)}\right)dy dx$} \nonumber \\
& \quad -\gamma \left\|\tilde{\mathbf{X}}(t)\right\|_{\mathbb{R}^n}^2 - 
2\gamma\mathbf{L}^T\mathbf{P}\tilde{\mathbf{X}}(t)\int\limits_{y_1}^{y_2}g(y)\beta(t,0,y)dy.
\end{align}
Denoting
\begin{subequations}
\begin{align}
M_{G_1} & = \max_{x,y \in [0,1]} |G_1(x,y)|, & M_{PL} & = \|\mathbf{PL}\|_{\mathbb{R}^n}^2, \\
M_{G_2} & = \max_{x,y \in [0,1]} |G_2(x,y)|, & M_g & = \|g\|^2_{L^2([y_1,y_2]; \mathbb{R})},
\end{align}
\end{subequations}
we can estimate using Cauchy-Schwarz and Young's inequalities
\begin{subequations}
\begin{align}
2\gamma\mathbf{L}^T\mathbf{P}\tilde{\mathbf{X}}(t)\int\limits_{y_1}^{y_2}g(y)\beta(t,0,y)dy & 
\leq \nonumber \\
\frac{\gamma}{2} \left\|\tilde{\mathbf{X}}(t)\right\|_{\mathbb{R}^n}^2 + 
2\gamma M_{PL}M_g\int\limits_0^1\beta(t,0,y)dy,  \\
\resizebox{.9\columnwidth}{!}{$\displaystyle 
\int\limits_0^1\int\limits_0^1 
\frac{2p_1(y)}{\lambda(x,y)}\alpha(t,x,y)G_1(x,y)\beta(t,0,y)\exp\left(-\delta\int\limits_0^x\frac{ds 
}{\lambda(s,y)}\right)dy dx$} & \leq \nonumber \\
\frac{\delta}{2}\int\limits_0^1\int\limits_0^1\frac{p_1(y)}{\lambda(x,y)}\alpha^2(t,x,y)\exp\left(-\delta\int\limits_0^x\frac{ds
 }{\lambda(s,y)}\right)dydx \nonumber \\
+ \int\limits_0^1 \frac{2}{\delta}\frac{p_1(y)}{m_\lambda}M_{G_1}^2 \beta^2(t,0,y)dy,
\end{align}
\begin{align}
\resizebox{.9\columnwidth}{!}{$\displaystyle 
\int\limits_0^1\int\limits_0^1 
\frac{2}{\mu(x,y)}\beta(t,x,y)G_2(x,y)\beta(t,0,y)\exp\left(\delta\int\limits_0^x\frac{ds}{\mu(x,y)}\right)dy
 dx$} & \leq \nonumber \\
\frac{\delta}{2}\int\limits_0^1\int\limits_0^1\frac{1}{\mu(x,y)}\beta^2(t,x,y)\exp\left(\delta\int\limits_0^x
\frac{ds}{\mu(x,y)}\right)dydx \nonumber \\
+ \int\limits_0^1 
\frac{2}{\delta}\frac{1}{m_\mu}M_{G_2}^2\exp\left(\frac{\delta}{m_\mu}\right)
\beta^2(t,0,y)dy,
\end{align}
\end{subequations}
which then gives
\begin{align}
\label{eq:lyapdest}
\dot{V}_\mathrm{tot}(t) & \leq -\frac{\delta}{2}V(t) -\int\limits_0^1 
\left(p_1(y)\exp\left(-\delta\int\limits_0^1\frac{
ds}{\lambda(s,y)}\right) \right. \nonumber \\
& \left.\qquad - 
R^2(y)\exp\left(\delta\int\limits_0^1\frac{ds}{\mu(s,y)}\right)\right)\alpha^2(t,1,y)dy 
\nonumber \\
& \quad -\int\limits_0^1 \left(1 - p_1(y)\frac{2M_{G_1}^2}{m_\lambda\delta} - 
\frac{2M_{G_2}^2}{m_\mu\delta}\exp\left(\frac{\delta}{m_\mu}\right) 
\right. 
\nonumber \\
& \left. \qquad - p_1(y)Q^2(y) - 2\gamma M_{PL}M_g
\right)\beta^2(t,0,y)dy 
\nonumber \\
& \quad  - \frac{\gamma}{2}\left\|\tilde{\mathbf{X}}(t)\right\|_{\mathbb{R}^n}^2.
\end{align}
In order to cancel out the term multiplying $\alpha^2(t,1,y)$ in \eqref{eq:lyapdest}, we choose
\begin{equation}
p_1(y)  = R^2(y)\exp\left(\delta\int\limits_0^1\left(\frac{1}{\lambda(s,y)} + 
\frac{1}{\mu(s,y)}\right)ds\right),
\end{equation}
inserting which to the term multiplying $\beta^2(t,0,y)$ in \eqref{eq:lyapdest} gives that 
$\dot{V}_\mathrm{tot}$ is negative definite if, for all $y \in [0,1]$,
\begin{align}
\resizebox{.9\columnwidth}{!}{$
\left(R(y)Q(y)\right)^2\exp\left(\delta\int\limits_0^1\left(\frac{1}{\lambda(s,y)}
 + \frac{1}{\mu(s,y)}\right)ds\right) +
\frac{2M_{G_2}^2}{\delta m_\mu}\exp\left(\frac{\delta}{m_\mu}\right)$}
\nonumber \\
\resizebox{.9\columnwidth}{!}{$+ 2\gamma M_{PL}M_g + R^2(y)\frac{2M_{G_1}^2}{\delta 
m_\lambda}\exp\left(\delta\int\limits_0^1\left(\frac{1}{\lambda(s,y)}
 + \frac{1}{\mu(s,y)}\right)ds\right)$} & \leq 1.
\label{eq:sstabcond}
\end{align}
Now, fixing $\delta > 0$ sufficiently small such that \eqref{eq:deltacond} holds and having $M_G$ 
satisfy \eqref{eq:MGcond}, we can fix $\gamma > 0$ sufficiently small such that 
\eqref{eq:sstabcond} is satisfied. Finally, as $R(y) \neq 0$ for all $y \in 
[0,1]$, all weights in the Lyapunov functional \eqref{eq:Vtot} are positive. Hence, 
$V_\mathrm{tot}$ is equivalent to the norm on $\mathbb{R}^n \times E_c^2$, and exponential 
stability of \eqref{eq:obstsmod}, \eqref{eq:obstsbc} on $\mathbb{R}^n \times E_c^2$ follows by 
negative-definiteness of $\dot{V}_\mathrm{tot}$.
\end{proof}
\end{lemma}

\begin{theorem}
\label{thm:cobs}
Under the conditions of Lemma~\ref{lem:bswp}, if
\eqref{eq:MGcond} and \eqref{eq:deltacond} hold, then the estimation error system 
\eqref{eq:lerr}, \eqref{eq:lerrbc} is exponentially stable on $\mathbb{R}^n \times E_c^2$.

\begin{proof}
The proof follows directly by Lemmas~\ref{lem:bswp} and \ref{lem:stab}, considering the  
backstepping transformations \eqref{eq:obsbs} that (for each fixed $y\in[0,1]$) are Volterra 
transformations of second 
kind (in $x$) for $(\alpha, \beta)$, and hence, boundedly invertible on $E_c^2$ by \cite[Thm 
2.3.5]{HocBook}.
\end{proof}
\end{theorem}

\subsection{Relation to Sylvester Equation-Based Observer Design} \label{obs:syl}

The purpose of this subsection is to show that, if the PDE part of \eqref{eq:infart}, 
\eqref{eq:infartbc} is stable, the backstepping approach adopted here to 
estimate the ODE-PDE cascade coincides with the Sylvester equation-based approach from 
\cite[Thm 4.4]{Nat21}. This is because the kernel equations \eqref{eq:kg1}, \eqref{eq:kg2} 
with boundary conditions \eqref{eq:kg1bc}, \eqref{eq:kg2bc} coincide with the Sylvester equation 
\cite[(4.4)]{Nat21} given by
\begin{equation}
\label{eq:sylv}
\Pi\mathbf{A} = A\Pi + B\mathbf{C},
\end{equation}
on $\mathbb{R}^n$, where $A, B$ are operators associated to the abstract Cauchy problem 
formulation of the PDE part of \eqref{eq:infart}, \eqref{eq:infartbc}. Moreover, this provides an 
alternative 
representation for $\pmb{\gamma}^1, \pmb{\gamma}^2$, as highlighted in the following 
proposition.
\begin{proposition}
\label{prop:sylv}
The gains $\pmb{\gamma}_1, \pmb{\gamma}_2$ appearing in \eqref{eq:P12} are connected to the 
solution of the Sylvester equation \eqref{eq:sylv} as
\begin{equation}
\begin{bmatrix}
\pmb{\gamma}^1 \\ \pmb{\gamma}^2
\end{bmatrix} = \Pi,
\end{equation}
and they have the alternative (equivalent) representations\footnote{Formulas 
\eqref{eq:gammas}--\eqref{eq:vk} are valid for eigenvalues of geometric multiplicity one. For 
exposition clarity and for not distracting the reader with algebraic details from the main result, 
which is representation of the backstepping kernels via solution to a Sylvester equation, this is 
(only) tacitly assumed in Proposition~\ref{prop:sylv}. In the 
general case, the solution to \eqref{eq:sylv} (in abstract form) is given in \cite[(4.6)]{Nat21}.}
\begin{subequations}
\label{eq:gammas}
\begin{align}
\pmb{\gamma}_1(x,y) & = \sum_{k=1}^n u^k(x,y)\mathbf{v}_k^*, \\
\pmb{\gamma}_2(x,y) & = \sum_{k=1}^n v^k(x,y)\mathbf{v}_k^*,
\end{align}
\end{subequations}
where 
\begin{equation}
\label{eq:ukvk}
\resizebox{.98\columnwidth}{!}{$\displaystyle
\begin{bmatrix}
u^k(x,y) \\ v^k(x,y)
\end{bmatrix} = \exp\left(\int\limits_0^x \begin{bmatrix}
\frac{-s_k}{\lambda(\xi,y)} & \frac{W(\xi,y)}{\lambda(\xi,y)} \\
-\frac{\theta(\xi,y)}{\mu(\xi,y)} & \frac{s_k}{\mu(\xi,y)}
\end{bmatrix}d\xi\right)
\begin{bmatrix}
Qv^k(0,y) \\ v^k(0,y)
\end{bmatrix}$},
\end{equation}
with
\begin{equation}
\label{eq:vk}
\resizebox{.98\columnwidth}{!}{$\displaystyle 
v^k(0,y) = \frac{F(y)\mathbf{Cv}_k}{\begin{bmatrix} -R(y) & 1
\end{bmatrix}\exp\left(\int\limits_0^1 \begin{bmatrix}
\frac{-s_k}{\lambda(\xi,y)} & \frac{W(\xi,y)}{\lambda(\xi,y)} \\
-\frac{\theta(\xi,y)}{\mu(\xi,y)} & \frac{s_k}{\mu(\xi,y)}
\end{bmatrix}d\xi\right)\begin{bmatrix}
Q(y) \\ 1
\end{bmatrix}}$},
\end{equation}
and $\mathbf{v}_k, s_k$ are the eigenvectors and the respective eigenvalues of 
$\mathbf{A}$. Furthermore, under Assumption~\ref{ass:obs},
\eqref{eq:gammas}--\eqref{eq:vk} 
are well-defined and the pair $\left(\mathbf{A}, 
\int\limits_{y_1}^{y_2}g(y)\pmb{\gamma}^2(0,y)dy\right)$ is detectable.

\begin{proof}
Applying \eqref{eq:sylv} to an eigenvector $\mathbf{v}_k$ of $\mathbf{A}$ and denoting the 
respective eigenvalue by $s_k$, we get 
\begin{equation}
\label{eq:sylveig}
\Pi s_k\mathbf{v}_k = A\Pi\mathbf{v}_k + B\mathbf{C}\mathbf{v}_k,
\end{equation}
the solution to which is
\begin{equation}
\label{eq:sylveigsol}
\Pi\mathbf{v}_k = (s_k - A)^{-1}B\mathbf{Cv}_k,
\end{equation}
where $(s_k - A)^{-1}B\mathbf{Cv}_k := (u^k,v^k)$ is the solution to the ODE (see 
\cite[Rem. 
10.1.5]{TucWeiBook})
\begin{subequations}
  \label{eq:sylode}
\begin{align}
s_k u^k(x,y) & = -\lambda(x,y)u^k_x(x,y) + W(x,y)v^k(x,y), \\
s_kv^k(x,y) & = \mu(x,y)v^k_x(x,y) + \theta(x,y)u^k(x,y),
\end{align}
\end{subequations}
with boundary conditions
\begin{subequations}
  \label{eq:sylodebc}
\begin{align}
u^k(0,y) & = Q(y)v^k(0,y), \\
v^k(1,y) & = R(y)u^k(1,y) + F(y)\mathbf{Cv}_k,
\end{align}
\end{subequations}
which, by denoting $u^k(x,y) = \pmb{\gamma}^1(x,y)\mathbf{v}_k$ and $v^k(x,y) = 
\pmb{\gamma}^2(x,y)\mathbf{v}_k$, coincides with \eqref{eq:kg1}, \eqref{eq:kg2}, 
\eqref{eq:kg1bc}, \eqref{eq:kg2bc}  applied to 
an eigenvector $\mathbf{v}_k$ of $\mathbf{A}$ with $s_k$ being the corresponding eigenvalue. 
The solution to \eqref{eq:sylode}, \eqref{eq:sylodebc} is then given by \eqref{eq:ukvk},
and the boundary conditions \eqref{eq:sylodebc} give
\begin{align}
\label{eq:FyCv}
F(y)\mathbf{Cv}_k & = \begin{bmatrix} -R(y) & 1
\end{bmatrix}\exp\left(\int\limits_0^1 \begin{bmatrix}
\frac{-s_k}{\lambda(\xi,y)} & \frac{W(\xi,y)}{\lambda(\xi,y)} \\
-\frac{\theta(\xi,y)}{\mu(\xi,y)} & \frac{s_k}{\mu(\xi,y)}
\end{bmatrix}d\xi\right) \nonumber \\
& \quad \times \begin{bmatrix}
Q(y) \\ 1
\end{bmatrix}v^k(0,y),
\end{align}
which gives the expression \eqref{eq:vk} for $v^k(0,y)$. In order to show that $v^k(0,y)$ is 
well-defined, we note that
\begin{equation}
G(s_k)\mathbf{Cv}_k = \int\limits_{y_1}^{y_2} g(y)v^k(0,y)dy,
\end{equation}
where $G(s_k)$ denotes the transfer function of the continuum $2\times 2$ system from 
$F(y)U(s) = v(s,1,y) - 
R(y)u(s,1,y)$ \\ to $Y(s) = \int\limits_{y_1}^{y_2} g(y)v(s,0,y)dy$ for $s = s_k$, which exists by 
the first part of Assumption~\ref{ass:obs}, and hence, $v^k(0,y)$ is well-defined. Moreover, by 
\cite[Prop. 4.7]{Nat21}, the necessary and sufficient condition for $\left(\mathbf{A}, 
\int\limits_{y_1}^{y_2}g(y)\pmb{\gamma}^2(0,y)dy\right)$ to be detectable is that 
$G(s_k)\mathbf{Cv}_k \neq 0$ for all (unstable) eigenvalues $s_k$ and respective eigenvectors 
$\mathbf{v}_k$ of $\mathbf{A}$, which is guaranteed by the second part of 
Assumption~\ref{ass:obs}.
\end{proof}
\end{proposition}

\section{Continuum-Based Observer Design for ODE - Large-Scale Hyperbolic PDE 
Cascades} \label{sec:art}

In this section, we present the class of large-scale collections of $2\times 2$ hyperbolic PDE 
systems considered and 
introduce a constructive approach to approximate them by a respective continuum. 
Thereafter, we utilize the continuum system constructed and apply the observer design from 
Section~\ref{sec:obs} to (approximately) estimate the state of a class of ODE - large-scale 
$2\times 2$ hyperbolic PDE systems by the continuum observer.

\subsection{Large-Scale Collections of $2\times 2$ Hyperbolic Systems} \label{art:sys}

Consider a large-scale collection of $2\times 2$ hyperbolic systems for $i \in \{1,\ldots,m\}$  of the 
form
\begin{subequations}
\label{eq:arti}
\begin{align}
u_t^i(t,x) + \lambda^i(x)u_x^i(t,x) & = w^i(x)v^i(t,x), \\
v_t^i(t,x) - \mu^i(x)v_x^i(t,x) & = \theta^i(x)u^i(t,x),
\end{align}
\end{subequations}
with boundary conditions
\begin{subequations}
\label{eq:artibc}
\begin{align}
u^i(t,0) & = q^iv^i(t,0), \\
v^i(t,1) & = r^iu^i(t,1) + f^i\mathbf{C}\mathbf{X}(t),
\end{align}
\end{subequations}
where the dynamics of $\mathbf{X}$ are driven by the ODE
\begin{equation}
\label{eq:ode}
\dot{\mathbf{X}}(t)  = \mathbf{A}\mathbf{X}(t),
\end{equation}
as in the ODE - continuum-PDE cascade \eqref{eq:infart}, \eqref{eq:infartbc}.
The parameters of \eqref{eq:arti}, \eqref{eq:artibc} satisfy the 
following assumption.
\begin{assumption}
\label{ass:m2x2}
The parameters of \eqref{eq:arti}, \eqref{eq:artibc} are such that $\lambda^i,\mu^i \in C^1([0,1]; 
\mathbb{R}), \theta^i,w^i \in C([0,1]; \mathbb{R})$, and $r^i,q^i,f^i$ $\in \mathbb{R}$, for all $i = 
1,\ldots,m$, where $\lambda^i(x), \mu^i(x) > 0$ for all $x \in [0,1]$.
\end{assumption}

As an output of \eqref{eq:arti}, \eqref{eq:artibc}, we consider a weighted average of the 
boundary values $v^i(t,0)$ of the form
\begin{equation}
\label{eq:artiout}
Y(t)=\frac{1}{m}\sum_{i=m_1}^{m_2} g^i v^i(t,0),
\end{equation}
for some weights $g^i$ and integers $1 \leq m_1 \leq 
m_2 \leq m$. This is motivated by non-invasive average measurements in the arterial network 
application considered in Section~\ref{sec:artappl}.

The PDE part of system \eqref{eq:arti}, \eqref{eq:artibc} can be written compactly as an $m+m$ 
system with 
states 
$\mathbf{u} = (u^1,\ldots,u^m), \mathbf{v} = (v^1,\ldots,v^m)$ as 
\begin{subequations}
\label{eq:nmm2}%
\begin{align}
	\mathbf{u}_t(t,x) + \pmb{\Lambda}(x)\mathbf{u}_x(t,x)   & 
	= 
	\mathbf{W}(x)\mathbf{v}(t,x), \\
	\mathbf{v}_t(t,x) - \mathbf{M}(x)\mathbf{v}_x(t,x)
	& = 
	\pmb{\Theta}(x)\mathbf{u}(t,x),
\end{align}
\end{subequations}
with boundary conditions 
\begin{subequations}
\label{eq:nmm2bc}
\begin{align}
\mathbf{u}(t,0) & = \mathbf{Q}\mathbf{v}(t,0), \\
\mathbf{v}(t,1) & =\mathbf{R}\mathbf{u}(t,1) + \mathbf{F}\mathbf{C}\mathbf{X}(t),
\end{align}
\end{subequations}
where $ \mathbf{F} = \left(f^i\right)_{i=1}^m$, while all the parameter matrices are diagonal and 
contain the respective parameters from 
\eqref{eq:arti}, \eqref{eq:artibc}, i.e.,
\begin{subequations}
\begin{align}
 \pmb{\Lambda}  &
= \operatorname{diag}(\lambda^1,\ldots,\lambda^m), & 
 \mathbf{M}  &
= \operatorname{diag}(\mu^1,\ldots,\mu^m), \\
 \mathbf{W}  &
= \operatorname{diag}(w^1,\ldots,w^m), & 
 \pmb{\Theta}  &
= \operatorname{diag}(\theta^1,\ldots,\theta^m), \\
 \mathbf{Q}  &
= \operatorname{diag}(q^1,\ldots,q^m), & 
 \mathbf{R}  &
= \operatorname{diag}(r^1,\ldots,r^m).
\end{align}
\end{subequations}

\subsection{Construction of Continuum Approximation} \label{art:cont}

In order to compare the solutions of a large-scale $2\times 2$ system with a continuum $2\times 
2$ system, we 
interpret the parameters of the large-scale $2\times 2$ systems as step functions in $y$ as 
follows (respectively for the other parameters)
\begin{equation}
\label{eq:lamstep}
\lambda^s(x,y) = \begin{cases}
\lambda^1(x), & y \in [0, 1/m], \\
\lambda^i(x), & y \in ((i-1)/m), i/m], \; i = 2,\ldots,m.
\end{cases}
\end{equation}
With the step function interpretation of the large-scale parameters, the accuracy of a continuum 
approximation can be measured based on, e.g., the smallest $\varepsilon > 0$ such that 
\begin{subequations}
\label{eq:infappr}
\begin{align}
\sup_{x\in[0,1]}\|\lambda^s(x,\cdot) - \lambda(x,y\cdot)\|_{L^2([0,1]; \mathbb{R})} \qquad  
\nonumber \\
\qquad + \sup_{x\in[0,1]}\|\lambda_x^s(x,\cdot) - \lambda_x(x,\cdot)\|_{L^2([0,1]; \mathbb{R})} & 
< 
\varepsilon, \label{eq:infappr1}
\\
\sup_{x\in[0,1]}\|\mu^s(x,\cdot) - \mu(x,\cdot)\|_{L^2([0,1]; \mathbb{R})} \qquad  \nonumber \\
+ \sup_{x\in[0,1]}\|\mu_x^s(x,\cdot) - \mu_x(x,\cdot)\|_{L^2([0,1]; \mathbb{R})} & < 
\varepsilon,
\\
\sup_{x\in[0,1]}\|w^{s}(x,\cdot) - W(x,\cdot)\|_{L^2([0,1]; \mathbb{R})} & < \varepsilon,  \\
\sup_{x\in[0,1]}\|\theta^{s}(x,\cdot) - \theta(x,\cdot)\|_{L^2([0,1]; \mathbb{R})} & < \varepsilon,  \\
\|q^{s}  - Q\|_{L^2([0,1]; \mathbb{R})} & < \varepsilon, \\
\|r^{s}  - R\|_{L^2([0,1]; \mathbb{R})} & < \varepsilon, \\
\|f^{s}  - F\|_{L^2([0,1]; \mathbb{R})} & < \varepsilon, \\
\|g^{s} - g\|_{L^2([y_1,y_2]; \mathbb{R})} & < \varepsilon.
\end{align}
\end{subequations}
Conversely, \eqref{eq:infappr} can be viewed as a guideline for constructing the continuum 
parameters, e.g., find $\lambda$ such that \eqref{eq:infappr1} holds for a given $\varepsilon$. 
Thus, considering a model function $f_\lambda$ such that $\lambda(x,y) = 
f_\lambda(x,y,\mathbf{c}_\lambda)$, where $\mathbf{c}_\lambda$ is a vector of adjustable 
parameters, the objective would be to find $\mathbf{c}_\lambda$ such that 
$f_\lambda(x,y,\mathbf{c}_\lambda)$ minimizes the left-hand side of \eqref{eq:infappr1}.

In practice, however, the $\sup$-norm is not the most convenient norm for solving optimization 
problems, so we rather consider approximation in the $L^2$-norm also in $x$. Hence, instead of 
the left-hand side of, e.g., \eqref{eq:infappr1}, we consider
\begin{equation}
\label{eq:lqrL2}
\|\lambda^s - \lambda\|_{E_c}^2 + \|\lambda^s_x - \lambda_x\|_{E_c}^2.
\end{equation}
For computational simplicity, we consider $\lambda$ to be a linear combination of tunable 
parameters $\mathbf{c}_\lambda$, i.e., 
\begin{equation}
\label{eq:lampar}
\lambda(x,y)  = \sum_k c_\lambda^k\varphi^k(x,y),
\end{equation}
for chosen functions $\varphi^k$. Inserting \eqref{eq:lampar} to \eqref{eq:lqrL2} and taking the 
gradient with respect $\mathbf{c}_\lambda$, we get
\begin{align}
\frac{\partial}{\partial \mathbf{c}_\lambda}\left(\left\|\lambda^s - \sum_k 
c_\lambda^k\varphi_k\right\|_{E_c}^2 + \left\|\lambda^s_x - \sum_k 
c_\lambda^k\varphi^k_x\right\|_{E_c}^2\right) & = \nonumber \\
-2\left(\mathbf{V}^\lambda - \mathbf{E}\mathbf{c}_\lambda + \mathbf{D}^\lambda - 
\mathbf{G}\mathbf{c}_\lambda\right),
\end{align}
where the elements of $\mathbf{V}^\lambda, \mathbf{E}, \mathbf{D}^\lambda, \mathbf{G}$ are 
given by
\begin{subequations}
\label{eq:L2lq}
\begin{align}
\mathbf{V}_k^\lambda & = \langle \varphi^k, \lambda^s \rangle_{E_c},  & 
\mathbf{E}_{k,\ell} & = \langle \varphi^k, \varphi^\ell \rangle_{E_c}, \\
\mathbf{D}_k^\lambda & = \langle \varphi_x^k, \lambda_x^s \rangle_{E_c}, &
\mathbf{G}_{k,\ell} & = \langle \varphi_x^k, \varphi_x^\ell \rangle_{E_c}.
\end{align}
\end{subequations}
The minimum of the quadratic function \eqref{eq:lqrL2} is obtained when its gradient is zero, 
which then gives
\begin{equation}
\label{eq:clsol}
\mathbf{c}_\lambda = \left(\mathbf{E} + \mathbf{G}\right)^{-1}\left(\mathbf{V}^\lambda + 
\mathbf{D}^\lambda\right).
\end{equation}

In order to expedite computations, we take $\varphi^k$ as orthogonal functions so that 
$\mathbf{E}$ is diagonal and $\mathbf{G}$ may also be sparse. For example, with Legendre 
polynomials, $\mathbf{G}$ can be written as band matrix with bandwidth $2$. Thus, a polynomial 
approximation of total order $M$ and potentially reduced-order $M_y$ in $y$ can be written as
\begin{equation}
\label{eq:lamleg}
\lambda(x,y) = \sum_{j=0}^{M_y}\sum_{i=0}^{M-j}c_\lambda^{(i,j)}L_i(x)L_j(y),
\end{equation}
where $L_i$ denotes Legendre polynomial of order $i$, and $k$ in 
\eqref{eq:lampar} can be connected to $(i,j)$ through
\begin{equation}
k(i,j) = i + Mj - \frac{(j-1)(j-2)}{2} + 2,
\end{equation}
where $k  = 1, \ldots,M(M_y+1) - \frac{M_y(M_y-1)}{2} + 1$. Now, in terms of finding the simplest 
continuum approximation, one can consider finding the smallest $M_y$ (given some $M$) for 
which the least-squares error \eqref{eq:lqrL2} is smaller than a chosen $\bar{\varepsilon} > 0$ 
(note that the $L^2$-norm only provides a lower bound for the $L^\infty$-norm in 
\eqref{eq:infappr}). Alternatively, 
one can consider a multi-objective optimization problem, where the goal is to minimize both the 
least-squares error \eqref{eq:lqrL2} and the approximation order $M_y$. Finally, an 
additional 
benefit of Legendre polynomials is that, if the parameter data are (piecewise) smooth, the 
polynomial approximation converges pointwise (see, e.g., \cite{BabHak19}). Hence, in such cases, 
the least-squares fit in the $L^2$ norm also leads to small approximation errors in 
\eqref{eq:infappr}.

Finding the polynomial fit for $\mu^s$ is analogous to $\lambda^s$. For $w^s$ and $\theta^s$, 
the process is similar, except that the $\partial_x$ terms are not present in the computations. For 
$q^s,r^s$, and $f^s$, the process is analogous to $\theta^s$ and $w^s$, except that the inner 
products in \eqref{eq:L2lq} are in 1-D instead of 2-D. Finally, we note that 
under Assumption~\ref{ass:loc}, the constructed $\lambda$ and $\mu$ have to be positive. If this 
is not satisfied by the least-squares coefficients \eqref{eq:clsol}, one option is to solve instead the 
constrained problem 
\begin{align}
\label{eq:clconstr}
& \min_{\mathbf{c}_\lambda} \left\|\left(\mathbf{E}+\mathbf{G}\right)\mathbf{c}_\lambda - 
\left(\mathbf{V}^\lambda+\mathbf{D}^\lambda\right)\right\|^2, \nonumber \\
\mathrm{s.t.} & \min_{(x,y) \in 
[0,1]^2}\sum_{j=0}^{M_y}\sum_{i=0}^{M-j}c_\lambda^{(i,j)}L_i(x)L_j(y) > 0,
\end{align}
where the unconstrained solution \eqref{eq:clsol} may provide a good initial guess. However, 
satisfying the constraint in \eqref{eq:clconstr} is by no means trivial, as evaluating the left-hand 
side already requires solving a constrained minimization problem. Regardless, considering that 
$\lambda^s > 0$, finding a sufficiently accurate fit $\lambda$ through \eqref{eq:clsol} may 
already satisfy $\lambda > 0$.

\begin{figure}[!htb]
\begin{center}
\includegraphics[width=\columnwidth]{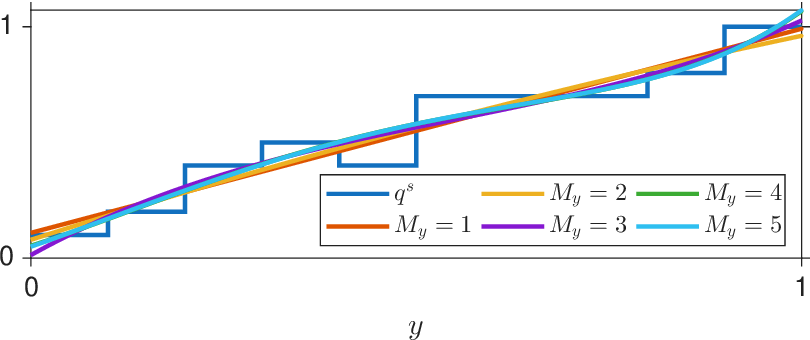}
\end{center}
\caption{The step function $q^s$ and polynomial fits of order $M_y = 1,\ldots,5$ for the $q^i$
data from \eqref{eq:qdata}.}
\label{fig:qtest}
\end{figure}

Based on the above, a systematic method for constructing continuum 
approximations for large-scale $2\times 2$ systems is presented in Algorithm~\ref{alg:contappr}. 
We note that it is left to the user's discretion to determine appropriate approximation accuracy, as 
the attainable level of accuracy is dependent on the large-scale parameters. In practice, one 
should simply compute 
approximations of different orders, and then choose the lowest-order approximation that gives 
acceptable approximation accuracy. For example, considering $q^i$ data given by
\begin{equation}
\label{eq:qdata}
q^i = [1, 2,4,5,4,7,7,7,8,10]/10,
\end{equation}
the step function $q^s$ according to \eqref{eq:qdata} and polynomial fits of order $M_y = 
1,\ldots,5$ are shown in Figure~\ref{fig:qtest}, where one can see that the polynomials of different 
orders do not differ much. The respective $L^2$ residuals are given in 
Table~\ref{tab:qtest}, where, analogously to Figure~\ref{fig:qtest}, the residuals do not differ 
much 
either. We note that, by construction, the $L^2$ residual decreases as $M_y$ increases, but the 
additional accuracy comes at the cost of increased complexity. Regardless, based on 
Table~\ref{tab:qtest} (and Figure~\ref{fig:qtest}), there is no apparent benefit 
of choosing a continuum approximation of order higher than $M_y=1$ in this case.

\begin{table}[!htb]
	\begin{center}
\resizebox{.99\columnwidth}{!}{
		\begin{tabular}{c | c c c c c}  
			$M_y$ & 1 & 2 & 3 & 4 & 5 \\
			\hline \\ [-6pt]
			$\|q^s-q\|_{L^2}$ & 0.0113  & 0.0110 & 0.0097 & 0.0094 & 0.0094
		\end{tabular}}
	\end{center}
	\caption{$L^2$ residuals for the polynomial fits in Figure~\ref{fig:qtest}.}
	\label{tab:qtest}
\end{table}

\begin{algorithm}[!htb]
	\DontPrintSemicolon
	\KwData{Parameters $\lambda^i,\mu^i,\theta^i,w^i,q^i,r^i,f^i,g^i,m_1,m_2$ of 
	\eqref{eq:arti}, \eqref{eq:artibc}, \eqref{eq:artiout}.}
	\KwResult{Continuum parameters $\lambda,\mu,\theta,W,Q,R,F,g$.}
	\Init{}{
		Construct step functions $\lambda^s,\mu^s,w^s,\theta^s,q^s,r^s,f^s,g^s$ based on the 
		parameter data as in \eqref{eq:lamstep}.\;
		Set $y_1 = (m_1-1)/m$ and $y_2 = m_2/m$. \;
		Choose approximation orders $M$ and $M_y$.
	}
	\While{continuum approximation not satisfactory}{
		Compute the matrices $\mathbf{E}, \mathbf{G}$, the vectors $\mathbf{D}^\lambda, 
		\mathbf{D}^\mu$, and $\mathbf{V}^{\lambda,\mu,w,\theta,q,r,f,g}$ as in 
		\eqref{eq:L2lq}. \;
		Compute tunable parameters as in \eqref{eq:clsol} and construct continuum parameters as in 
		\eqref{eq:lamleg}. \;
		Check approximation accuracy \eqref{eq:infappr} and least-squares errors as in 
		\eqref{eq:lqrL2}. \;
		If necessary, adjust $M_y$ (and $M$), otherwise deem continuum approximation satisfactory.
	}
	\If{$\lambda > 0$ (resp. $\mu > 0$) violated}{
		Solve $\mathbf{c}_\lambda$ (resp. $\mathbf{c}_\mu$) from \eqref{eq:clconstr} and insert to 
		\eqref{eq:lamleg}.
	}
	\label{alg:contappr}
	\caption{Construction of continuum approximation for a large-scale $2\times 2$ system 
	\eqref{eq:arti}, \eqref{eq:artibc}, \eqref{eq:artiout}.}
\end{algorithm}

\subsection{Observer Design via Continuum Approximation}

We provide in the following theorem the theoretical guarantees of estimation of the state of the 
large-scale system \eqref{eq:arti}--\eqref{eq:ode} with measurement \eqref{eq:artiout}, by the 
continuum observer \eqref{eq:lobs}, \eqref{eq:lobsbc} with measurement \eqref{eq:artiout}.

\begin{theorem}
\label{thm:caobs}
Consider a system \eqref{eq:arti}--\eqref{eq:ode} where the PDE part is exponentially stable and 
the ODE part is marginally stable, with parameters satisfying 
Assumption~\ref{ass:m2x2}, and output \eqref{eq:artiout}.
Construct a respective continuum observer \eqref{eq:lobs}, \eqref{eq:lobsbc} (with measurement 
obtained from \eqref{eq:artiout}) such that the 
continuum PDE parameters satisfy Assumptions~\ref{ass:loc} and \ref{ass:obs}, and 
\eqref{eq:infappr}. Then, for all $t \geq 0$,
\begin{align}
\label{eq:cobserr}
\left\|  \begin{bmatrix}
\tilde{\mathbf{X}}(t) \\ \tilde{u}(t) \\ \tilde{v}(t)
\end{bmatrix} \right\|_{\mathbb{R}^n \times E_c^2}
  & \leq
M_e e^{-\omega_e t}\left\| 
\begin{bmatrix}
\tilde{\mathbf{X}}_0 \\ \tilde{u}_0 \\ \tilde{v}_0
\end{bmatrix} \right\|_{\mathbb{R}^n\times E} \nonumber \\
  & \qquad + \delta\left\| 
\begin{bmatrix}
\mathbf{X}_0 \\ \mathbf{u}_0 \\ \mathbf{v}_0
\end{bmatrix} \right\|_{\mathbb{R}^n\times E},
\end{align}
for some $M_e,\omega_e > 0$, where $\delta$ depends continuously on $\varepsilon > 
0$ in \eqref{eq:infappr} such that $\delta \to 0$ as
$\varepsilon\to 0$, and we denote $\tilde{\mathbf{X}} = \mathbf{X} -
\hat{\mathbf{X}}, \tilde{u} = \mathcal{F}_m\mathbf{u} - \hat{u},
\tilde{v} = \mathcal{F}_m\mathbf{v} - \hat{v}$.

\begin{proof}
Denoting $z^m = (\mathbf{X}, \mathcal{F}_m\mathbf{u}, \mathcal{F}_m\mathbf{v})^T$, there 
exist operators $\mathbb{T}_t^m, \Psi_t^m$ such that the solution to 
\eqref{eq:arti}--\eqref{eq:artode} (with the $\mathbb{R}^m$ PDE states transformed to $E_c$ via 
the isometric transform $\mathcal{F}_m$) with output \eqref{eq:artiout} can be written as
\begin{align}
\label{eq:zm}
z^m(t) & = \mathbb{T}_t^mz_0^m, &
Y(t) &  = \Psi_t^m z_0^m.
\end{align}
Respectively, consider a virtual continuum system of the form \eqref{eq:infart}--\eqref{eq:out} 
based on the continuum parameters constructed in the statement of the theorem. There exist 
operators $\mathbb{T}_t, \Psi_t$ such that the solution to and output of  
\eqref{eq:infart}--\eqref{eq:out} can be written as
\begin{align}
  \label{eq:virtc}
z^c(t) & = \mathbb{T}_t z_0^m, &
Y^c(t) &  = \Psi_t z_0^m,
\end{align}
where the initial condition of the virtual continuum system has been
chosen as $z_0^m$. Moreover, introducing the input-to-state map
$\Phi_t$ corresponding to the output injection operator $[\mathbf{L},
P_1, P_2]^T$ in \eqref{eq:lobs}, the solution to the cascade system of
\eqref{eq:zm} and  \eqref{eq:lobs}, \eqref{eq:lobsbc} (with measurement \eqref{eq:artiout})
can be written as 
\begin{equation}
  \label{eq:pocas}
\begin{bmatrix}
z^m(t) \\ \hat{z}(t)
\end{bmatrix} = 
\begin{bmatrix}
\mathbb{T}_t^m & 0 \\ -\Phi_t\Psi_t^m & \mathbb{T}_t + \Phi_t\Psi_t
\end{bmatrix} 
\begin{bmatrix}
z_0^m \\ \hat{z}_0
\end{bmatrix}.
\end{equation}
Respectively, the solution to the virtual cascade system of  \eqref{eq:virtc}
and \eqref{eq:lobs}, \eqref{eq:lobsbc} (with measurement \eqref{eq:out}) can be written as 
\begin{equation}
  \label{eq:pocasc}
\begin{bmatrix}
z^c(t) \\ \hat{z}^c(t)
\end{bmatrix} = 
\begin{bmatrix}
\mathbb{T}_t & 0 \\ -\Phi_t\Psi_t & \mathbb{T}_t + \Phi_t\Psi_t
\end{bmatrix} 
\begin{bmatrix}
z_0^m \\ \hat{z}_0
\end{bmatrix}.
\end{equation}
Now, writing 
\begin{equation}
\tilde{z} = \hat{z} - z^m = (\hat{z}^c - z^c) - (z^m - z^c) - (\hat{z}^c - \hat{z}),
\end{equation}
we have by Theorem~\ref{thm:cobs} and \eqref{eq:pocas},
\eqref{eq:pocasc} that 
\begin{align}
  \|\tilde{z}(t)\|_{\mathbb{R}^m \times E_c^2}
  & \leq M_e e^{-\omega_e t}\tilde{z}_0 + \|(\mathbb{T}^m_t
    - \mathbb{T}_t)z_0^m\|_{\mathbb{R}^n\times E_c^2} \nonumber \\
  & \qquad + \|\Phi_t(\Psi_t^m
    - \Psi_t)z_0^m\|_{\mathbb{R}^n\times E_c^2},                                                                  
\end{align}
where the terms depending on $z_0^m$ are uniformly bounded in time by analogous/dual 
arguments to \cite[(70)--(73)]{HumBek26c}, and they tend to zero as $\varepsilon \to 0$ by 
analogous arguments to \cite[(74)]{HumBek26c} (that rely on the proof of \cite[Thm 
2.5]{HumBek26c}). Hence, the estimate 
\eqref{eq:cobserr} follows.
\end{proof}
\end{theorem}

We note that the result of Theorem~\ref{thm:caobs} is qualitative in the sense that, in 
\eqref{eq:cobserr}, the continuous dependence  of $\delta$ on $\varepsilon$ from 
\eqref{eq:infappr} 
is not explicitly provided. In order to quantify $\delta$ in terms of $\varepsilon$, one could employ 
similar tools as, e.g., in \cite[Prop. 1]{HumBek26b} to explicitly estimate the continuum 
approximation error. However, the estimate obtained in \cite[Prop. 1]{HumBek26b} seems to be 
highly conservative, so that, practically, a more efficient manner to quantify the achievable 
accuracy of 
the estimation error is via numerical investigations, as we do for the numerical example in 
Section~\ref{num:ex1}.

\section{Numerical Implementation and Simulations} \label{sec:num}

\subsection{Spectral Approximation of the Continuum PDE} \label{num:spec}

Consider the PDE part of \eqref{eq:infart}  in weak form with a test
function $\phi \in E_c$ such that $\phi \in C^\infty([0,1]; \mathbb{R})$ and $\phi(0), \phi(1) 
\neq 0$ (with $E_{c}$ inner 
products)
\begin{subequations}
\label{eq:infartwf}  
\begin{align}
\left\langle u_t,\phi  \right\rangle + \left\langle
  \lambda u_x, \phi  \right\rangle   & = \left\langle Wv,\phi \right\rangle,  \\
\left\langle v_t,\phi \right\rangle - \left\langle \mu v_x, \phi  \right\rangle
& = \left\langle \theta u, \phi  \right\rangle.
\end{align}
\end{subequations}
Integrating by parts in $\left\langle \lambda u_x,\phi \right\rangle$
(respectively in $\left\langle \mu v_x,\phi \right\rangle$ gives
\begin{align}
\left\langle \lambda u_x,\phi \right\rangle
& = \int\limits_0^1 \left([\lambda u\phi](1,y) -
     [\lambda u\phi](0,y) \right)dy - \left\langle u,(\lambda \phi)_x
     \right\rangle \nonumber \\
& = \int\limits_0^{1}\left( [\lambda u\phi](1,y) - Q(y)[\lambda
                                  v\phi](0,y) \right)dy \nonumber \\
  & \qquad- \left\langle u,\lambda_x\phi \right\rangle - \left\langle
    u,\lambda\phi_x \right\rangle, 
\end{align}
where we also used the boundary conditions \eqref{eq:infartbc}. Now,
writing a spectral approximation for $u$ and $v$ as 
\begin{subequations}
\begin{align}
u(t,x,y) & \approx \sum_{j=0}^{N_y}\sum_{i=0}^{N_x}
           c_{i,j}^u(t)\phi_i(x)\phi_j(y), \\
v(t,x,y) &  \approx \sum_{j=0}^{N_y}\sum_{i=0}^{N_x}
           c_{i,j}^v(t)\phi_i(x)\phi_j(y),
\end{align}
\end{subequations}
for some $N_y, N_x \in \mathbb{N}$, where $\phi_i, \phi_j$ are
orthogonal basis functions, the weak form \eqref{eq:infartwf} for
arbitrary $\phi_{k,\ell} = \phi_k\phi_\ell$ becomes 
\begin{subequations}
\label{eq:wffull}  
\begin{align}
\dot{c}^u_{k,\ell}(t) \left\langle \phi_{k,\ell},\phi_{k,\ell} \right\rangle +
 \sum_{j=0}^{N_y}\sum_{i=0}^{N_x}
           c_{i,j}^u(t) \int\limits_0^1
  [\lambda\phi_{k,\ell}\phi_{i,j}](1,y)dy \nonumber \\
  - \sum_{j=0}^{N_y}\sum_{i=0}^{N_x}
           c_{i,j}^v(t) \int\limits_0^1
  Q(y)[\lambda \phi_{k,\ell}\phi_{i,j}](0,y)dy \nonumber \\
  -  \sum_{j=0}^{N_y}\sum_{i=0}^{N_x}
           c_{i,j}^u(t) \left\langle
  \phi_{i,j},\lambda_{x}\phi_{k,\ell}+\lambda\partial_x\phi_{k,\ell}
  \right\rangle  & = \nonumber \\
  \sum_{j=0}^{N_y}\sum_{i=0}^{N_x}
           c_{i,j}^v(t)\left\langle \phi_{i,j}, W\phi_{k,\ell}
  \right\rangle, \\
  \dot{c}^v_{k,\ell}(t) \left\langle \phi_{k,\ell},\phi_{k,\ell}
  \right\rangle  + \sum_{j=0}^{N_y}\sum_{i=0}^{N_x}
           c_{i,j}^v(t) \int\limits_0^1
  [\mu\phi_{k,\ell}\phi_{i,j}](0,y)dy \nonumber \\
  - \sum_{j=0}^{N_y}\sum_{i=0}^{N_x}
           c_{i,j}^u(t) \int\limits_0^1
  R(y)[\mu\phi_{k,\ell}\phi_{i,j}](1,y)dy
  \nonumber \\
  -  \int\limits_0^1
  F(y)[\mu\phi_{k,\ell}](1,y)dy \mathbf{CX}(t)
  \nonumber \\
  +  \sum_{j=0}^{N_y}\sum_{i=0}^{N_x}
           c_{i,j}^v(t) \left\langle
  \phi_{i,j},\mu_{x}\phi_{k,\ell}+\mu\partial_x\phi_{k,\ell}
  \right\rangle  & = \nonumber \\
  \sum_{j=0}^{N_y}\sum_{i=0}^{N_x}
           c_{i,j}^u(t)\left\langle \phi_{i,j}, \theta\phi_{k,\ell}
  \right\rangle.
\end{align}
\end{subequations}
Now, introducing an indexing function $(i,j)=\mathcal{I}(k)$ for $k = 
1,2,\ldots,(N_x+1)(N_y+1)$, e.g., as
\begin{equation}
\label{eq:indfun}
(i,j) = \mathcal{I}(k) = \left(\operatorname{mod}(k-1,N_x+1), \left\lfloor \frac{k-1}{N_x+1} 
\right\rfloor\right),
\end{equation}
we can write the coefficients $c_{i,j}^u, c_{i,j}^v$ as vectors 
\begin{equation}
\mathbf{c}_u = 
\left(c^u_{\mathcal{I}(k)}\right)_{k=1}^{(N_x+1)(N_y+1)}, \quad \mathbf{c}_v = 
\left(c^v_{\mathcal{I}(k)}\right)_{k=1}^{(N_x+1)(N_y+1)}.
\end{equation}
Finally, testing \eqref{eq:wffull} over all
$\phi_{\ell,k}$, we get an ODE 
\begin{equation}
  \label{eq:infartODE}  
\begin{bmatrix}
\mathbf{M} & 0 \\ 0 & \mathbf{M}
\end{bmatrix} 
\begin{bmatrix}
\dot{\mathbf{c}}_{u} \\ \dot{\mathbf{c}}_v
\end{bmatrix} = 
\begin{bmatrix}
\mathbf{K}^{uu} & \mathbf{K}^{vu} \\ \mathbf{K}^{uv} & \mathbf{K}^{vv}
\end{bmatrix} 
\begin{bmatrix}
\mathbf{c}_u \\ \mathbf{c}_v
\end{bmatrix} + 
\begin{bmatrix}
0 \\ \mathbf{B}
\end{bmatrix} \mathbf{CX}(t),
\end{equation}
with parameters
\begin{subequations}
\begin{align}
\mathbf{M} & = \operatorname{diag}\left(\left\langle \phi_{\mathcal{I}(k)},\phi_{\mathcal{I}(k)}
             \right\rangle\right)_{k=1}^{(N_x+1)(N_y+1)}, \\
\mathbf{K}^{uu}_{k,\ell} & = \left\langle
  \phi_{\mathcal{I}(\ell)},\lambda_{x}\phi_{\mathcal{I}(k)}+\lambda\partial_x\phi_{\mathcal{I}(k)}
  \right\rangle \nonumber \\
& \qquad  - \int\limits_0^1
  [\lambda\phi_{\mathcal{I}(k)}\phi_{\mathcal{I}(\ell)}](1,y)dy, \\
\mathbf{K}^{vu}_{k,\ell} & = \left\langle 
  \phi_{\mathcal{I}(\ell)}, W\phi_{\mathcal{I}(k)}
  \right\rangle + \int\limits_0^1
  Q(y)[\lambda \phi_{\mathcal{I}(k)}\phi_{\mathcal{I}(\ell)}](0,y)dy, \\
\mathbf{K}^{uv}_{k,\ell} & = \left\langle \phi_{\mathcal{I}(\ell)}, \theta\phi_{\mathcal{I}(k)}
  \right\rangle + \int\limits_0^1
  R(y)[\mu\phi_{\mathcal{I}(k)}\phi_{\mathcal{I}(\ell)}](1,y)dy, \\
\mathbf{K}^{vv}_{k,\ell} & = -\left\langle
  \phi_{\mathcal{I}(\ell)},\mu_{x}\phi_{\mathcal{I}(k)}+\mu\partial_x\phi_{\mathcal{I}(k)}
  \right\rangle \nonumber \\
& \qquad  - \int\limits_0^1
  [\mu\phi_{\mathcal{I}(k)}\phi_{\mathcal{I}(\ell)}](0,y)dy, \\
\mathbf{B}_k & = \int\limits_0^1
  F(y)[\mu\phi_{\mathcal{I}(k)}](1,y)dy,
\end{align}
\end{subequations}
for $k,\ell = 1,\ldots,(N_x+1)(N_y+1)$. We take the orthogonal test/basis functions $\phi_{i,j}$ 
as Legendre polynomials, i.e., $\phi_{i,j}(x,y) = L_i(x)L_j(y)$, which satisfy the conditions 
stated in the beginning of this subsection.

For implementing the continuum observer \eqref{eq:lobs}, \eqref{eq:lobsbc}, the ODE resulting 
from the spectral approximation is analogous to \eqref{eq:infartODE}, with the addition of the 
output injection term
\begin{equation}
\begin{bmatrix}
\mathbf{P}^1 \\ \mathbf{P}^2
\end{bmatrix}\begin{bmatrix}
0 & \mathbf{C}^v
\end{bmatrix}\begin{bmatrix}
\tilde{\mathbf{c}}_u \\ \tilde{\mathbf{c}}_v
\end{bmatrix},
\end{equation}
where
\begin{subequations}
\begin{align}
\mathbf{P}^1_k & = \langle P_1, \phi_{\mathcal{I}(k)}\rangle, \qquad
\mathbf{P}^2_k = \langle P_2, \phi_{\mathcal{I}(k)}\rangle, \\
\mathbf{C}^v_k & = \int\limits_{y_1}^{y_2}g(y)\phi_{\mathcal{I}(k)}(0,y)dy.
\end{align}
\end{subequations}
Thus, a spectral approximation for the PDE part of the observer is given by
\begin{align}
  \label{eq:infartobsODE}  
\begin{bmatrix}
\mathbf{M} & 0 \\ 0 & \mathbf{M}
\end{bmatrix} 
\begin{bmatrix}
\dot{\hat{\mathbf{c}}}_{u} \\ \dot{\hat{\mathbf{c}}}_v
\end{bmatrix} & = 
\begin{bmatrix}
\mathbf{K}^{uu} & \mathbf{K}^{vu} + \mathbf{P}^1\mathbf{C}^v \\ \mathbf{K}^{uv} & 
\mathbf{K}^{vv} + \mathbf{P}^2\mathbf{C}^v
\end{bmatrix} 
\begin{bmatrix}
\hat{\mathbf{c}}_u \\ \hat{\mathbf{c}}_v
\end{bmatrix} \nonumber \\
& \qquad + 
\begin{bmatrix}
0 \\ \mathbf{B}
\end{bmatrix} \mathbf{C}\hat{\mathbf{X}}(t) 
- \begin{bmatrix}
0 & \mathbf{P}^1\mathbf{C}^v \\ 0 & \mathbf{P}^2\mathbf{C}^v
\end{bmatrix} \begin{bmatrix}
\mathbf{c}_u \\ \mathbf{c}_v
\end{bmatrix}.
\end{align}

\subsection{Academic Numerical Example} \label{num:ex1}

For an academic example, consider a large-scale $2\times 2$ system of the form \eqref{eq:arti}, 
\eqref{eq:artibc} with parameters, for $i = 
1,\ldots,10$ and $x \in [0,1]$,
\begin{subequations}
\label{eq:exparam}
\begin{align}
\lambda^i(x) & = 1, & \mu^i(x) & = 1, \\
w^i(x) & =  c_w^ix(x+1), & \theta^i(x) & = c_\theta^ix ,
\end{align}
\end{subequations}
with $f^i = 1$ for $i = 1,\ldots,10$, and
\begin{subequations}
\label{eq:exparam2}
\begin{align}
q^i & = [1, 2,4,5,4,7,7,7,8,10]/10, \\
r^i & = [9,8,6,4,4,3,2,2,-1,0]/10, \\
c_w^i & = [-3, -3, -2, -2, 1, 0, 1, 2, 3, 6]/20, \\
c_\theta^i & = [12, 11, 14,13,15,16,18,19,18,20]/20.
\end{align}
\end{subequations}
In order to derive a continuum counterpart of the form \eqref{eq:infloc1}, \eqref{eq:infloc2}, 
\eqref{eq:infartbc} we construct continuum approximations for the parameters 
\eqref{eq:exparam}, \eqref{eq:exparam2}, employing Algorithm~\ref{alg:contappr} from
Section~\ref{art:cont}. Choosing $M=3$ and $M_y = 1$,\footnote{As $\theta^i$ and $w^i$ are 
separable in $x$ and $i$, it would be enough to find 1-D polynomial fits to the $c_\theta^i$ and 
$c_w^i$ data, respectively, but we carry out the process in 2-D regardless.} we get the  
following continuum parameters, for $x,y \in [0,1]$
\begin{subequations}
\label{eq:exparamc}
\begin{align}
\lambda(x,y) & = 1, \qquad  \mu(x,y) = 1, \\
W(x,y) & =  \left(\frac{906}{2000}y - \frac{423}{2000} \right)x(x+1), \\
\theta(x,y) & = \left(\frac{486}{1000}y + \frac{537}{1000}\right)x, \\
Q(y) & = \frac{882}{1000}y + \frac{109}{1000}, \\
 R(y)&  = -\frac{1026}{1000}y + \frac{883}{1000}, \qquad F(y) = 1. \label{eq:paramcR}
\end{align}
\end{subequations}
The ODE-part is taken as a harmonic signal generator approximating the proximal flow 
waveform based on the Fourier coefficients listed in \cite[Table 6.1]{StePhD} scaled by a factor 
of $10^3$ (see Section~\ref{sec:artappl} for details). The 
gain $\mathbf{L}$ is computed using LQR (\texttt{lqr} in MATLAB) with state weight $10^4$ and 
input weight $1$. Moreover, $\pmb{\gamma}_1, \pmb{\gamma}_2$ are computed by solving 
\eqref{eq:sylode}, \eqref{eq:sylodebc} using power series of order $40$ in $x$\footnote{Solving 
\eqref{eq:sylode}, \eqref{eq:sylodebc} accurately for the higher 
frequencies appears to require high approximation order in $x$.} and reduced order 
$4$ in $y$. For the measurement \eqref{eq:artiout}, we take $m_1 = 1$ and $m_2 = m$ with $g^i 
= 1$ for all $i \in \{1,\ldots,m\}$, so that in the respective continuum measurement \eqref{eq:out} 
we have $y_1 = 0$ and $y_2 = 1$ with $g(y) = 1$ for all $y \in [0,1]$. We note that $R$ in 
\eqref{eq:paramcR} violates the condition $R(y) \neq 0$ for all $y \in [0,1]$, so that 
Lemma~\ref{lem:stab} cannot be used to check stability of the PDE part of the system. Instead, 
stability of the continuum PDE for parameters \eqref{eq:exparamc} was validated in a numerical 
simulation.

\begin{figure}[!htb]
\begin{center}
\includegraphics[width=\columnwidth]{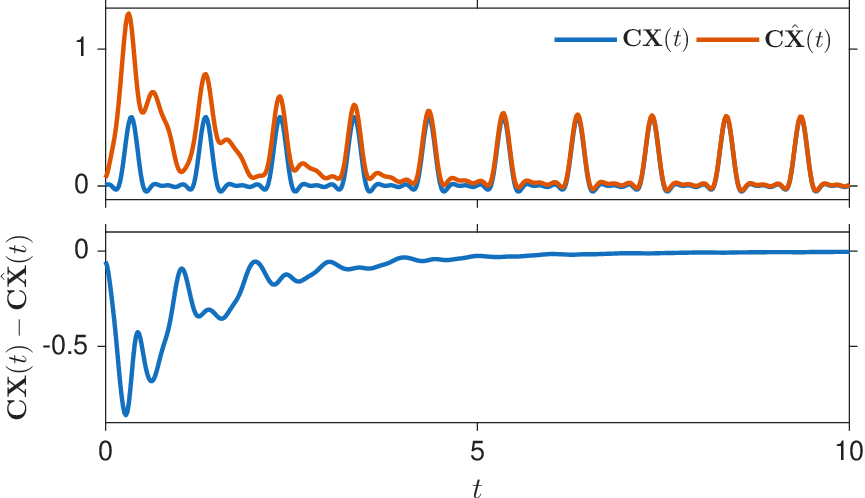}
\end{center}
\caption{Comparison of $\mathbf{CX}$ and the estimate $\mathbf{C}\hat{\mathbf{X}}$ for $N_y = 
M_y = 1$.}
\label{fig:ex11}
\end{figure}

For the simulation, the PDE parts of the plant \eqref{eq:arti}--\eqref{eq:ode}
and observer \eqref{eq:lobs}, \eqref{eq:lobsbc} are implemented based on the spectral 
approximation method of Section~\ref{num:spec}, where $N_x = 14$ is taken for both systems 
and $N_y = 1$ is taken 
for the continuum observer (as the continuum parameters \eqref{eq:exparamc} are first-order 
polynomials in $y$). All plant channels are initialized to $u_0^i(x) = v_0^i(x) = \sin(\pi x)$ while 
the PDE part of the observer is initialized to zero and the ODE part to one. The simulation results 
are shown in Figure~\ref{fig:ex11}, where the 
ODE output $\mathbf{CX}(t)$ and its estimate are shown for $t \in [0,10]$, along with the 
difference of the two. It can be seen that, after the transient errors, the estimation error tends 
close to zero. However, the estimation error does not converge exactly to zero due to the 
continuum observer, where continuum approximation errors persist between the plant and the 
observer.

\begin{figure}[!htb]
\begin{center}
\includegraphics[width=\columnwidth]{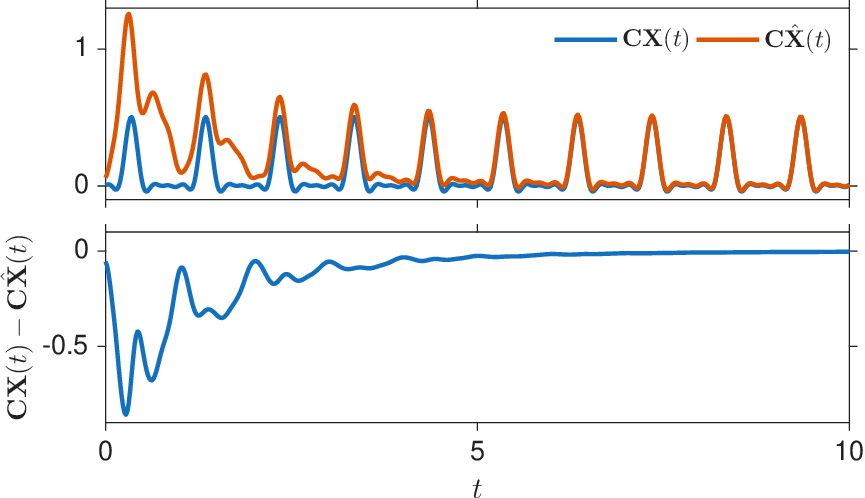}
\end{center}
\caption{Comparison of $\mathbf{CX}$ and the estimate $\mathbf{C}\hat{\mathbf{X}}$ for $N_y = 
M_y = 3$.}
\label{fig:ex12}
\end{figure}

\begin{figure}[!htb]
\begin{center}
\includegraphics[width=\columnwidth]{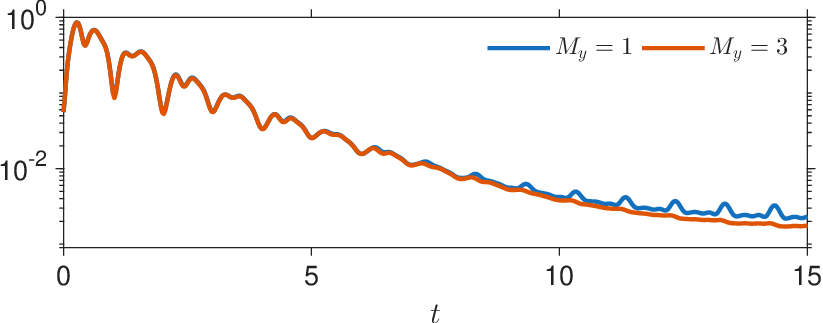}
\end{center}
\caption{Comparison of $|\mathbf{CX} - \mathbf{C}\hat{\mathbf{X}}|$ for $M_y = N_y  = 1$ and 
$M_y = N_y  = 3$.}
\label{fig:ex1log}
\end{figure}

For comparison, we repeat the above process with $M_y = 3$, which, in order to retain accuracy, 
necessitates to increase the order of power series for solving $\pmb{\gamma}_1, 
\pmb{\gamma}_2$ from \eqref{eq:sylode}, \eqref{eq:sylodebc} to order $8$ in $y$ (and $40$ in 
$x$ as before). Moreover, the continuum observer is implemented with reduced order $N_y = 
M_y = 3$ (and $N_x = 14$ as before). The respective simulation results are shown in 
Figure~\ref{fig:ex12}. which seem virtually identical to those in Figure~\ref{fig:ex11}. An 
explanation for 
such similar results can be found in Table~\ref{tab:ex1}, where one can see that the least-squares 
fits of order $M_y = 1$ already have rather small approximation errors. While the 
higher-order fits for 
$M_y = 3$ do achieve smaller errors, the absolute improvement in the approximation 
accuracy is minor. Finally, closer inspection of the estimation errors $\mathbf{CX} - 
\mathbf{C}\hat{\mathbf{X}}$ in Figure~\ref{fig:ex1log} shows that the higher order approximation 
does lead to smaller continuum approximation errors in the simulation, but similarly to the 
residuals of the least-squares fits in Table~\ref{tab:ex1}, the absolute improvement in the 
approximation accuracy is very minor. Hence, one can argue that this minor improvement in the 
approximation accuracy is not worth the increase in computational complexity when $M_y = 3$ 
compared to that when $M_y = 1$, even if both are still considerably lower than the number 
$m=10$ of subsystems in the original, large-scale $2\times 2$ system. Note that a 
spectral-based implementation of an observer that is designed directly for the large-scale system 
\eqref{eq:arti}, \eqref{eq:artibc}, would result in significantly higher computational complexity as, 
in such a case, the value of $M_y$ (or rather $m$) would a priori be fixed to $10$. 

A comparison of the considered continuum observers and the exact, large-scale observer is 
provided in Table~\ref{tab:ex1comp}, which shows that the continuum observers for $M_y = 1$ 
and $M_y = 3$ reduce the number of unknowns in the spectral approximations of the PDE 
observer by $80\%$ and $60\%$, respectively, as compared with implementation of an observer 
directly for the  large-scale system with $m=10$. As a trade-off, the continuum observers 
result in small residuals in the estimation error $|\mathbf{CX} - \mathbf{C}\hat{\mathbf{X}}|$. 
Nevertheless, the residuals are considered minor, especially considering the case $M_y = 1$, 
compared to the decrease in 
computational complexity of the respective spectral approximations.
\begin{table}[!htb]
	\begin{center}
		\begin{tabular}{c | c c c} 
			& continuum & continuum & large-scale \\
			 & $M_y = 1$ & $M_y = 3$ & $m = 10$ \\ 
			\hline \\ [-10pt]
			\# states & 60 & 120 & 300 \\
			residual &  0.0022 & 0.0015 & 0
		\end{tabular}
	\end{center}
	\caption{Comparison of the number of states in the spectral approximation of the 
	continuum/large-scale PDE observers for $N_x = 14$ and the respective residuals of the 
	estimation error $|\mathbf{CX} - \mathbf{C}\hat{\mathbf{X}}|$.}
	\label{tab:ex1comp}
\end{table}

\begin{table}[!htb]
	\begin{center}
\resizebox{\columnwidth}{!}{
		\begin{tabular}{c | c c c c} 
			$M_y$ & $\displaystyle \sup_{x \in [0,1]}\|W(x,\cdot)-w^s(x,\cdot)\|_{L^2}$ & 
			$\displaystyle \sup_{x \in [0,1]}\|\theta(x,\cdot)-\theta^s(x,\cdot)\|_{L^2}$ & 
			$\|Q-q^s\|_{L^2}$ & 
			$\|R-r^s\|_{L^2}$ \\ 
			\hline \\ [-10pt]
			1 & 0.0077 & 0.0038 & 0.0113 & 0.0128  \\
			3 & 0.0052 & 0.0031 & 0.0097 & 0.0098 
		\end{tabular}}
	\end{center}
	\caption{Approximation errors for the polynomial parameters of orders $M_y = 1$ and 
	$M_y = 3$ in $y$.}
	\label{tab:ex1}
\end{table}

\section{Application to Blood Flow Arterial Networks} \label{sec:artappl}

\subsection{Control-Theoretic Model of an Arterial Network} \label{art:thm}

\begin{figure}[!htb]
\begin{center}
\begin{tikzpicture}
\node at (-1.25,5) {
\begin{tabular}{c}
Aortic flow and \\ pressure $(Q_\mathrm{a},P_\mathrm{a})$
\end{tabular}
};
\node at (3.5,7.5) {
\begin{tabular}{c}
	Arterial flow and \\ pressure $(Q_i,P_i)$
\end{tabular}
};
\draw [ultra thick] (-2.5,4.5) -- (0,4.5);
\draw [ultra thick, dashed,-stealth] (0,4.5) -- (2,4);
\draw[draw=black, ultra thick] (2,6) rectangle (5,7) node[pos=.5] {
\begin{tabular}{c}
		$2\times 2$ hyperbolic \\ system $i=1$
\end{tabular}
};
\draw [ultra thick, -stealth] (0,4.5) -- (2,6.5);
\draw[draw=black, ultra thick] (2,4.5) rectangle (5,5.5) node[pos=.5] {
\begin{tabular}{c}
		$2\times 2$ hyperbolic \\ system $i=2$
\end{tabular}
};
\draw [ultra thick, -stealth] (0,4.5) -- (2,5);
\node at (3.5,4.1) {$\pmb{\vdots}$};
\draw[draw=black, ultra thick] (2,2.5) rectangle (5,3.5) node[pos=.5] {
\begin{tabular}{c}
	$2\times 2$ hyperbolic \\ system $i=m$
\end{tabular}
};
\draw [ultra thick, -stealth] (0,4.5) -- (2,3);
\draw [thick, stealth-] (2,2) -- (5,2);
\node at (2,1.6) {$x=1$};
\node at (5,1.6) {$x=0$};
\end{tikzpicture}
\end{center}
\caption{Schematic of a parallel network model of pressure and	flow dynamics in the arterial tree 
based on \cite[Figure 1]{SwaXuD09}.}
\label{fig:afex}
\end{figure}
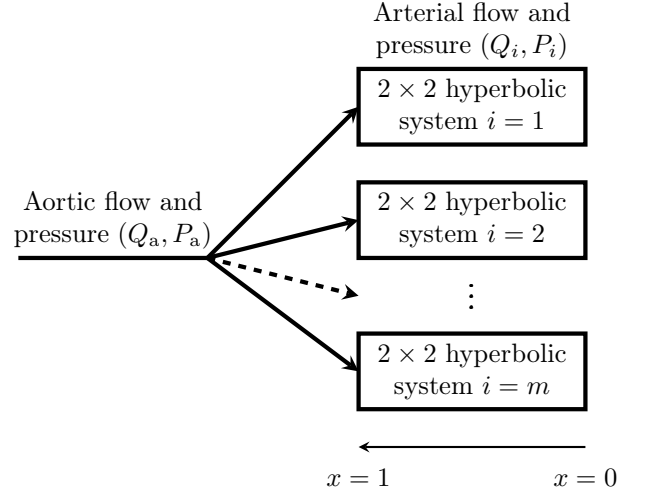

Figure~\ref{fig:afex} shows a network of $m$ peripheral arteries to which aortic blood flow enters. 
The main goal of such a network model (see \cite{SwaXuD09}) is to estimate aortic flow/pressure, 
which may be possible to measure only invasively, from peripheral measurements taken 
non-invasively (see also, e.g., \cite{BikPhd, GyuSot23}, for a description of practical ways of 
obtaining peripheral flow/pressure). 

A simplified, linearized model for blood flow in each 
peripheral artery, can be recast as a $2\times 2$ hyperbolic PDE system, which gives rise to a 
large-scale collection of $m$, $2 \times 2$ hyperbolic systems, as considered here. Therefore, a 
linearized model of 
blood flow in a single artery $i \in \{1,\ldots,m\}$  can be written as \eqref{eq:arti}, 
\eqref{eq:artibc} such that the parameters satisfy Assumption~\ref{ass:m2x2}. The aortic 
flow/pressure  waveform dynamics  are governed by \eqref{eq:ode} corresponding to, e.g., 
Fourier approximation of the waveform; see, e.g., \cite[Table 6.1]{StePhD} 
for blood flow waveform 
harmonics. We can consider
\begin{subequations}
\label{eq:artode}
\begin{align}
\mathbf{A} & = \operatorname{blkdiag}\left(0, \left(\begin{bmatrix}
0 & \omega_j \\ -\omega_j & 0
\end{bmatrix}
\right)_{j=1}^p\right), \\
\mathbf{C} & = \begin{bmatrix}
a_0 & b_1 & a_1 & \cdots & b_p & a_p
\end{bmatrix},
\end{align}
\end{subequations}
for some $(\omega_j)_{j=1}^p$, so by assigning the initial condition as $\mathbf{X}(0) = 
\mathbf{X}_0 = 
(1,0,1,\ldots,0,1)$, the output of the ODE is 
\begin{equation}
\mathbf{C}\mathbf{X}(t) = a_0 + \sum_{j=1}^p \left(a_j\cos(\omega_jt) + b_j\sin(\omega_jt)\right).
\end{equation}
A peripheral flow 
(resp. pressure) measurement 
from a single artery, or as an average of certain arteries, is possible to obtain non-invasively see, 
e.g., \cite{BikPhd, GyuSot23}. Thus, since the state $v$ consists of a linear combination of 
pressure and flow quantities (see the discussion below), it can be described in the form 
\eqref{eq:artiout}.

Following the control-theoretic, traffic 
flow-inspired model from \cite{Bek23} (that originates in 
\cite{StePhD, ReyMer09}), the nonlinear 
$2\times 2$ hyperbolic 
PDE for a single artery $i \in \{1,\ldots,m\}$ is of the form
\begin{subequations}
\label{eq:bf}
\begin{align}
A^i_t(x,t) & = A^i_x(x,t)V^i(x,t) + A^i(x,t)V^i_x(x,t), \\
V^i_t(x,t) & = V^i(x,t)V^i_x(x,t) + \frac{1}{\rho}P^i_x(A^i(x,t)) - K_r\frac{V^i(x,t)}{A^i(x,t)},
\end{align}
\end{subequations}
with boundary conditions
\begin{subequations}
\label{eq:bfbc}
\begin{align}
A^i(1,t)V^i(1,t) & = Q^i_{\mathrm{a}}(t), \label{eq:bfbc1} \\
P^i(A^i(0,t)) & = R_T^i A^i(0,t) V^i(0,t),
\end{align}
\end{subequations}
where $A^i > 0$ is the section area of the artery, $V^i > 0$ is the average blood
speed, $\rho > 0$ is the blood density, $K_r > 0$ is the friction parameter
related to blood viscosity, and $R_T^i 
\geq 0$ denotes total terminal resistance. The pressure function can be 
expressed as
\begin{equation}
\label{eq:PAi}
P^i(A^i) = \frac{\beta}{A_0^i}\left(\sqrt{A^i}-  \sqrt{A_0^i}\right),
\end{equation}
where $A_0^i$ is the reference arterial area at rest and $\beta = h^iE\sqrt{\pi}b$, where $h^i> 0$ 
is the artery wall thickness, $E > 0$ 
is Young's modulus, and $b$ is a positive parameter related to Poisson ratio of arterial 
wall. 
Then noting also that $Q^i_\mathrm{a}$, where $Q^i_\mathrm{a}$=$p_i Q_\mathrm{a}$ for some 
$p_i\in[0,1]$, with 
$\sum_{i=1}^m
p_i=1$ (see \cite{StePhD}), we get a system with input being $Q_\mathrm{a}$ only.

In order to obtain system \eqref{eq:arti}, \eqref{eq:artibc}, one has to transform system 
\eqref{eq:bf}, \eqref{eq:bfbc} to Riemann coordinates, as in \cite[Sect. III]{Bek23}, and to then 
linearize the system obtained as in  \cite[Sect. IV]{SinBek24} (using also \eqref{eq:PAi} to 
express the cross-sectional area $A$ as function of pressure $P$, as well as set $\mathbf{CX} 
= 
Q_{\rm a}$). Here, for simplicity, following 
also the modeling approach that employs time-invariant dynamics (in fact, transfer functions) 
from \cite{SwaXuD09} we consider linearization around a reference solution. The accuracy of this 
linearization, given the time-varying behavior of $Q_\mathrm{a}$, due to the nature of blood 
flows, has to be further explored. However, our modeling approach here serves as first step 
towards a more realistic model.

\subsection{Aortic Flow Estimation} \label{num:ex2}

Consider the blood flow arterial network model \eqref{eq:bf}--\eqref{eq:PAi} and linearize the 
system around the 
(constant) steady state $(V^i_*,A^i_*) = (0,A_0^i)$, where the artery is at rest. Towards that end, 
defining the Riemann variables \cite[(16), (17)]{Bek23}
\begin{subequations}
\label{eq:Riemann}
\begin{align}
u_R^i(A^i,V^i) & = V^i + 2\sqrt{\frac{2\beta}{\rho A_0^i}} \sqrt[4]{A^i}, \\
v_R^i(A^i,V^i) & = V^i - 2\sqrt{\frac{2\beta}{\rho A_0^i}} \sqrt[4]{A^i},
\end{align}
\end{subequations}
and linearizing the resulting system \cite[(20), (21)]{Bek23} leads to dynamics of the form 
\eqref{eq:arti}, where 
\begin{subequations}
\begin{align}
\lambda^i & = \frac{5u_*^i + 3v_*^i}{8},  \qquad \mu^i = -\frac{3u_*^i + 5v_*^i}{8}, \\
w^i(x) & = \tilde{w}^i\exp\left(\frac{\tilde{\sigma}^i}{\lambda^i}x +\frac{\tilde{\psi}^i}{\mu^i}x 
\right), \\
\theta^i(x) & = \tilde{\theta}^i\exp\left(-\frac{\tilde{\sigma}^i}{\lambda^i}x 
-\frac{\tilde{\psi}^i}{\mu^i}x \right),
\end{align}
\end{subequations}
where, according to \cite[(15)--(17)]{SinBek24}, we denote 
\begin{subequations}
\begin{align}
\tilde{\sigma}^i = \tilde{\theta}^i & = -\kappa^i\frac{3u_*^i + 5v_*^i}{(u_*^i - v_*^i)^5}, \\
\tilde{w}^i = \tilde{\psi}^i & = 
\kappa^i\frac{5u_*^i + 3v_*^i}{(u_*^i - v_*^i)^5},
\end{align}
\end{subequations}
with $\kappa^i = \frac{2^9K_r\beta^2}{\rho^2A_0^i{}^2}$, and $(u,v)$ are the transformed error 
variables
\begin{subequations}
\label{eq:trans}
\begin{align}
u^i(\cdot,x) & = \left(u_R^i - u_*^i\right)\exp\left(-\frac{\tilde{\sigma}^i}{\lambda^i}x\right), \\
v^i(\cdot,x) & = \left(v_R^i - v_*^i\right)\exp\left(\frac{\tilde{\psi}^i}{\mu^i}x\right),
\end{align}
\end{subequations}
around the steady state $(u_*^i, v_*^i) = (u_R^i(A_*^i,V_*^i),v_R^i(A_*^i,V_*^i))$. The 
parameters for the boundary conditions 
\eqref{eq:artibc} are obtained by linearizing \cite[(11)--(13)]{SinBek24} and by applying the 
transforms \eqref{eq:trans}.

For simulations, we consider ten arteries (i.e., $m=10$) where we vary the reference radius $r_0$ 
from $5.05$ 
mm to $5.5$ mm, so that the respective reference areas are $A_0^i = \pi r_0^i{}^2$ with $r_0 \in 
\{5.05,5.1,\ldots,5.5\}\cdot10^{-3}$. The other 
parameters are taken the same for all arteries based on the values used in \cite[Sect. 
III]{SinBek24}, i.e., $\rho = 1060$ kg/m$^3$, $K_r = 8\pi\mu/\rho$ with blood viscosity $\mu = 
0.0035$ Pa$\cdot$s, $h = 0.5$ mm, $E = 4\cdot 10^5$ N/m$^2$, $b = 4/3$, and $R_T 
=1.33\cdot 10^8$ N$\cdot$s/m$^5$. Similarly to Section~\ref{num:ex1}, the continuum 
parameters are 
constructed using  Algorithm~\ref{alg:contappr} for $M=3$ and $M_y = 3$, and the continuum 
observer is implemented with reduced order $N_y  =3$ in $y$, whereas $N_x = 14$ is taken in 
the $x$-direction for both the $m+m$ system and the continuum observer. The $m+m$ system is 
initialized as $u_0(x) = \frac{3}{4}\exp\left(\frac{\tilde{\sigma}^i}{\lambda^i}x\right)$ and $v_0(x) 
= 
-\frac{1}{2}\exp\left(-\frac{\tilde{\psi}^i}{\mu^i}x\right)$, and the ODE is the same as in 
Section~\ref{num:ex1}, i.e., taken of the form \eqref{eq:artode} based on the first four 
harmonics from \cite[Table 6.1]{StePhD} 
(with a slight adjustment in the zero-frequency harmonic to obtain positive flow values). 
Moreover, we take $Q_\mathrm{a}^i(t) = 0.1\cdot\mathbf{CX}(t)$ for all $i \in 
\{1,\ldots,10\}$, i.e., 
each considered artery receives ten percent of the aortic flow. For the output \eqref{eq:artiout}, 
we take $m_1 = 4$ and $m_2 = 7$ 
with $(g^i)_{i=4}^7 = (1,2,4,3)$, so that for $g$ in \eqref{eq:lobs} we  
have $y_1 = 0.3$ and $y_2 = 0.7$ with $g(y) \approx  -273.4y^3 +  375.0y^2 -155.9y +  
21.33$. This may correspond to a scenario in which the flow in peripheral arteries $6$ and $7$ 
mainly contribute to the average flow measurement. The 
PDE part of the observer 
is initialized to zero and the ODE part is initialized to 
one. 

\begin{figure}[!htb]
\begin{center}
\includegraphics[width=\columnwidth]{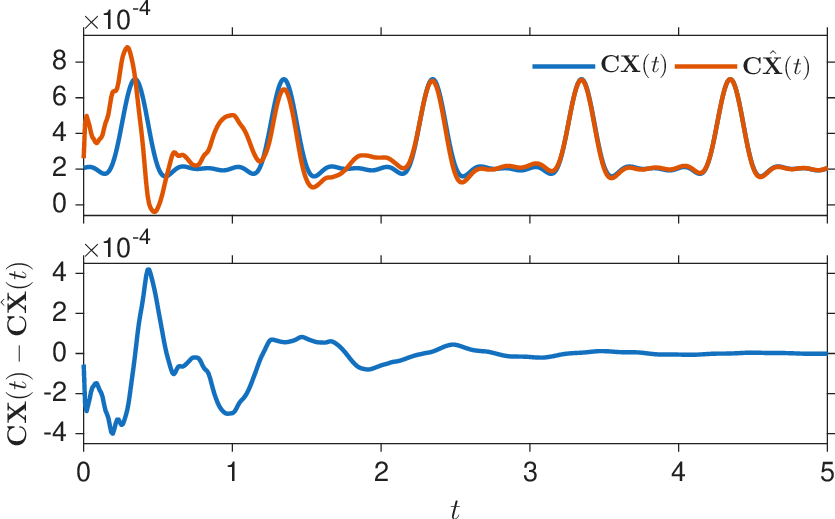}
\end{center}
\caption{Comparison of $Q_\mathrm{a} = \mathbf{CX}$ and the estimate 
$\mathbf{C}\hat{\mathbf{X}}$ for $M_y = N_y = 3$. The flow is given in  cubic meters per 
second.}
\label{fig:art}
\end{figure}

The simulation results are shown in Figure~\ref{fig:art}, where the ODE output and its 
estimate 
are shown for $t \in [0,5]$, along with the difference of the two. Compared to 
Section~\ref{num:ex1}, the 
estimation error decays much faster, which is due to the faster decay rate of the PDE dynamics 
here. As in Section~\ref{num:ex1}, the estimation error does not tend exactly to zero due to 
mismatches 
between the actual $m+m$ system and the continuum observer, albeit the estimation error 
caused by this mismatch is too small to be visible in Figure~\ref{fig:art}.

\subsection{Aortic Pressure Estimation}

In order to estimate aortic pressure, which may be an even more practically relevant problem, we 
modify the boundary condition \eqref{eq:bfbc1} to 
$P^i(A^i(1)) = P_\mathrm{a}$. In the Riemann variables 
\eqref{eq:Riemann} (using the inverse transformation \cite[(19)]{Bek23}), we have 
\begin{equation}
P^i(u^i,v^i) = \frac{\rho}{2^5}(u^i-v^i)^2 - \frac{\beta}{\sqrt{A_0^i}}.
\end{equation}
Linearization of $P^i$ around the same steady state as in the previous case gives rise to the 
boundary condition 
\begin{equation}
v^i(t,1) = r^iu^i(t,1) + f^iP_\mathrm{a}(t),
\end{equation}
where, denoting by $P_u$ and $P_v$ the partial derivatives of $P$ with respect to $u$ and $v$, 
respectively,
\begin{subequations}
\begin{align}
r^i & = -\frac{P^i_u(u_*^i,v_*^i)}{P^i_v(u_*^i,v_*^i)}\exp\left(\frac{\tilde{\sigma}^i}{\lambda^i} 
+\frac{\tilde{\psi}^i}{\mu^i} \right), \\
f^i & = \frac{1}{P^i_v(u_*^i,v_*^i)}\exp\left(\frac{\tilde{\psi}^i}{\mu^i} \right).
\end{align}
\end{subequations}
We set in all arteries $P_\mathrm{a}(t) = \mathbf{CX}(t)$, which is motivated by the assumption 
of pressure continuity at the bifurcation point as in \cite{SwaXuD09, StePhD}.
Otherwise the blood flow arterial network model is the same as in Section~\ref{num:ex2}, except 
the output, which is the same as in Section~\ref{num:ex1}. For the continuum approximation, we 
take $M = 3$ and $M_y = 1$ in Algorithm~\ref{alg:contappr} as well as $N_y=1$ for the numerical 
implementation of the observer ($N_x = 14$ is the same as previously). The
gain $\mathbf{L}$ is computed using LQR with state weight $100$ and input weight $1$, and the 
initial conditions 
for the $m+m$ system are taken as $4\exp\left(\frac{\tilde{\sigma}^i}{\lambda^i}x\right)$ and 
$v_0(x) = -\exp\left(-\frac{\tilde{\psi}^i}{\mu^i}x\right)$.

\begin{figure}[!htb]
\begin{center}
\includegraphics[width=\columnwidth]{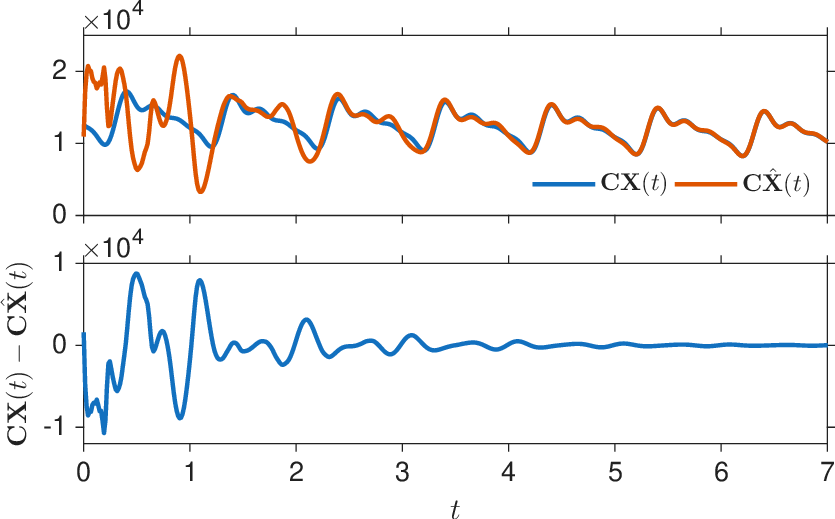}
\end{center}
\caption{Comparison of $P_\mathrm{a} = \mathbf{CX}$ and the estimate 
$\mathbf{C}\hat{\mathbf{X}}$ for $M_y = N_y = 1$. The pressure is given in Pascals.}
\label{fig:art2}
\end{figure}

The simulation results are shown in Figure~\ref{fig:art2}, where the reference aortic pressure 
$P_\mathrm{a}(t) = \mathbf{CX}(t)$ is designed to mimic the abdominal aorta waveform from 
\cite[Figure 6.14]{StePhD} by a (vanishing) Fourier series. The simulation results are consistent 
with 
those obtained in the previous examples.

\section{Conclusions and Future Work} \label{sec:conc}

We developed a continuum backstepping-based observer design for a class of ODE - 
continuum-PDE 
cascade systems and utilized it for (approximate) state 
estimation of a large-scale system counterpart, which originates from a linearized model of a 
blood flow arterial network. Moreover, we introduced  an approach to construct optimal 
continuum 
approximations for a given large-scale system, as well as a spectral-based approach to 
implement the resulting continuum observer. We provided theoretical, qualitative performance 
estimates for the estimation error systems and obtained consistent simulation results when the 
developed 
observer was applied to aortic flow/pressure estimation. In the present paper, the parameters of 
the arterial network model were assumed to be known. A next step, towards a more realistic 
modeling scenario as in \cite{SwaXuD09}, is to incorporate parameter adaptation into the state 
estimation, which may be pursued in future work.

\appendix

\setcounter{equation}{0}
\renewcommand{\theequation}{A.\arabic{equation}}

\section{Derivation of Kernel Equations} \label{app:ker}

Differentiating \eqref{eq:obsbs} with respect to $x$ and using the Leibniz rule, we get
\begin{subequations}
\begin{align}
u_x(t,x,y) & = \nonumber \\
 \alpha_x(t,x,y) - N^{11}(x,x,y)\alpha(t,x,y) \nonumber \\
- N^{12}(x,x,y)\beta(t,x,y)  -\int\limits_0^x N_x^{11}(x,\xi,y)\alpha(t,\xi,y)d\xi \nonumber \\
- \int\limits_0^x N_x^{12}(x,\xi,y)\beta(t,\xi,y)d\xi - 
\pmb{\gamma}_x^1(x,y)\tilde{\mathbf{X}}(t), \\
v_x(t,x,y) & = \nonumber \\
 \beta_x(t,x,y) - N^{21}(x,x,y)\alpha(t,x,y) \nonumber \\
- N^{22}(x,x,y)\beta(t,x,y) - \int\limits_0^xN_x^{21}(x,\xi,y)\alpha(t,\xi,y)d\xi \nonumber \\
- \int\limits_0^x N_x^{22}(x,\xi,y)\beta(t,\xi,y)d\xi - 
\pmb{\gamma}_x^2(x,y)\tilde{\mathbf{X}}(t).
\end{align}
\end{subequations}
From \eqref{eq:obsbs} and \eqref{eq:obstsmod}, we get
\begin{subequations}
\begin{align}
  \tilde{u}_t(t,x,y)
  & = \nonumber \\
-\lambda(x,y)\alpha_x(t,x,y) + G_1(x,y)\beta(t,0,y) \nonumber \\
+ \resizebox{.9\columnwidth}{!}{$\displaystyle  \int\limits_0^x
                 N^{11}(x,\xi,y)\left(\lambda(\xi,y)\alpha_\xi(t,\xi,y) - G_1(\xi,y)\beta(t,0,y)\right)d\xi$}
                 \nonumber \\
 -  
 \resizebox{.9\columnwidth}{!}{$\displaystyle\int\limits_0^xN^{12}(x,\xi,y)\left(\mu(\xi,y)\beta_\xi(t,\xi,y)
  + 
  G_2(\xi,y)\beta(t,0,y)\right)d\xi$} \nonumber \\
  -
    \pmb{\gamma}^1(x,y)\left(\left(\mathbf{A} + 
    \mathbf{L}\int\limits_{y_1}^{y_2}g(y)\pmb{\gamma}^2(0,y)dy 
\right)\tilde{\mathbf{X}}(t) \right. \nonumber \\ 
- \left.
\mathbf{L}\int\limits_{y_1}^{y_2}g(y)\beta(t,0,y)dy,\right),
\end{align}
\begin{align}
  \tilde{v}_t(t,x,y)
  & = 
\nonumber \\
\mu(x,y)\beta_x(t,x,y) + G_2(x,y)\beta(t,0,y) \nonumber \\
+
    \resizebox{.9\columnwidth}{!}{$\displaystyle 
    \int\limits_0^xN^{21}(x,\xi,y)\left(\lambda(\xi,y)\alpha_\xi(t,\xi,y) - 
    G_1(\xi,y)\beta(t,0,y)\right)d\xi $}
    \nonumber \\
 -  
 \resizebox{.9\columnwidth}{!}{$\displaystyle\int\limits_0^xN^{22}(x,\xi,y)\left(\mu(\xi,y)\beta_\xi(t,\xi,y)
  + 
 G_2(\xi,y)\beta(t,0,y)\right)d\xi$} \nonumber \\
- \pmb{\gamma}^2(x,y)\left(\left(\mathbf{A} +
\mathbf{L}\int\limits_{y_1}^{y_2}g(y)\pmb{\gamma}^2(0,y)dy 
\right)\tilde{\mathbf{X}}(t) \right. \nonumber \\ 
- \left.
\mathbf{L}\int\limits_{y_1}^{y_2}g(y)\beta(t,0,y)dy\right),
\end{align}
\end{subequations}
where (respectively for the other terms)
\begin{equation}
\resizebox{\columnwidth}{!}{$\begin{aligned}
\int\limits_0^x N^{11}(x,\xi,y)\lambda(\xi,y)\alpha_\xi(t,\xi,y)d\xi & = \nonumber \\
N^{11}(x,x,0)\lambda(x,y)\alpha(t,x,y) - N^{11}(x,0,y)\lambda(0,y)\alpha(t,0,y) \nonumber \\
-\int\limits_0^x \left(N^{11}_\xi(x,\xi,y)\lambda(\xi,y) + 
N^{11}(x,\xi,y)\lambda_\xi(\xi,y)\right)\alpha(t,\xi,y)d\xi.
\end{aligned}$}
\end{equation}
Thus, the kernel equations \eqref{eq:obsk}, \eqref{eq:obskbc} follow by setting \eqref{eq:lerr}, 
\eqref{eq:lerrbc} and \eqref{eq:obstsmod}, \eqref{eq:obstsbc} equal.

\end{document}